\documentclass[10pt]{scrartcl}
\usepackage{amsmath}
\usepackage{amssymb}
\usepackage{latexsym}
\usepackage{array}
\usepackage{yfonts}
\usepackage{fancyhdr}
\usepackage{epic, eepic}
\usepackage{graphicx,graphics}
\usepackage{color}
\usepackage{enumerate}
\DeclareMathAlphabet{\mathpzc}{OT1}{pzc}{m}{it}
\DeclareMathOperator{\vol}{vol}

\DeclareMathOperator{\Id}{Id}

\DeclareMathOperator{\diam}{diam}

\newtheorem{theorem}{Theorem}[section]
\newtheorem{defi}[theorem]{Definition}
\newtheorem{lemma}[theorem]{Lemma}

\newtheorem{cor}[theorem]{Corollary}
\newtheorem{rem}[theorem]{Remark}
\newtheorem{fig}{Figure}
\numberwithin{equation}{section}

\newcommand{\R}{\mathbb{R}}

\newcommand{\N}{\mathbb{N}}

\newcommand{\beg}{\begin{equation}}
\newcommand{\en}{\end{equation}}
\newcommand{\beqn}{\begin{eqnarray}}
\newcommand{\eeqn}{\end{eqnarray}}
\newcommand{\beqnn}{\begin{eqnarray*}}
\newcommand{\eeqnn}{\end{eqnarray*}}

\newcommand{\K}{\mathrm{k}}

\author{}
\title{}
\date{}
\begin{document}
\maketitle
\vspace{-3.1cm}
\begin{center}{\LARGE Immersions with bounded second fundamental form}\end{center}
\vspace{7mm}
\begin{center} Patrick Breuning \footnote{P.\ Breuning was supported by the DFG-Forschergruppe
\emph{Nonlinear Partial Differential Equations: Theoretical and Numerical Analysis}.
The contents of this paper were part of the author's dissertation, which was written at
Universit\"{a}t Freiburg, Germany.}
\\ Institut f\"{u}r Mathematik der Goethe Universit\"{a}t Frankfurt am Main \\ Robert-Mayer-Stra{\ss}e 10,
D-60325 Frankfurt am Main, Germany \\email: breuning@math.uni-frankfurt.de \end{center} \vspace{1mm}
\begin{abstract}
\begin{center} {\textbf{Abstract}} \end{center}
\noindent We first consider immersions on compact manifolds with uniform $L^p$-bounds on the second fundamental
form and uniformly bounded volume. We show compactness in
arbitrary dimension and codimension, generalizing a classical result of J.\ Langer. In the second part, this
result is used to deduce a localized version, being more convenient for many applications,
such as convergence proofs
for geometric flows.
\end{abstract}
\begin{section}{Introduction}
In \cite{langer} J.\ Langer investigated compactness of immersed surfaces
in $\R^3$ admitting uniform bounds on the second fundamental form and the area of the surfaces.
For a given sequence $f^i:\Sigma^i\rightarrow \R^3$, there exist after passing to a subsequence
a limit surface
$f:\Sigma\rightarrow \R^3$ and diffeomorphisms $\phi^i:\Sigma\rightarrow \Sigma^i$, such that
$f^i\circ\phi^i$ converges in the $C^1$-topology to $f$.
In particular, up to diffeomorphism, there are only finitely many
manifolds admitting such an immersion.
The finiteness of topological types was generalized by K.\ Corlette in \cite{corlette} to immersions of arbitrary
dimension and codimension. Moreover, the compactness theorem was generalized by S.\ Delladio in
\cite{delladio} to hypersurfaces of arbitrary dimension. \\ \\
The general case, that
is compactness in arbitrary dimension and codimension, is the first main theorem of this paper: \\
\begin{theorem}(Compactness theorem for immersions on compact manifolds) \label{compactness0} $ $\\ \\
Let $q$ be a point in $\R^n$, $m$ a positive integer, $p>m$,
and $\mathcal{A}, \mathcal{V}>0$ constants. Let $\mathfrak{F}$ be the set of
all mappings $f:M\rightarrow \R^n$ with the following properties:
\begin{itemize}
\item $M$ is an $m$-dimensional, compact manifold (without boundary)
\item $f$ is an immersion in $W^{2,p}(M,\R^n)$ with
\beqn \label{boundsecond} \|A(f)\|_{L^{p}(M)} &\leq& \mathcal{A} \hspace{2cm}\\
 \vol(M) &\leq& \mathcal{V}\\ q &\in& f(M).\label{qinfM} \eeqn
\end{itemize} Then for every sequence $f^i:M^i\!\rightarrow\!
\R^n$ in $\mathfrak{F}$ there exist a subsequence $f^j$, a mapping
$f:M\!\rightarrow\! \R^n$ in \nolinebreak$\mathfrak{F}$, and a sequence of
diffeomorphisms $\phi^j: M\!\rightarrow\! M^j$, such that $f^j \circ
\phi^j$ converges in the $C^1$-topology to $f$.
\end{theorem}
Here, the $L^p$-norm for the second fundamental form and the volume is measured
with respect to the volume measure induced by $f$.
Having shown Theorem \ref{compactness0}, we will use the
Nash embedding to generalize the result to complete Riemannian manifolds as target.
For a definition of the $C^1$-topology see \cite{hirsch}, p.\ 34--35.
The assumption $q\in f(M)$ ensures that the immersions $f^i$ do not diverge uniformly. This
can we weakened to $f(M)\cap K\neq \emptyset$ for a
fixed compact set $K\subset \R^n$.
In the case of an $L^\infty$-bound on the second fundamental form,
the assumption $\vol(M)\leq \mathcal{V}$ is equivalent to a bound on
the diameter $\diam(M)\leq \mathcal{D}$.
The theorem can easily be generalized to higher order convergence, provided we assume
uniform bounds $\|\nabla^kA\|_{L^\infty(M)}\leq \mathcal{A}_k$ for all covariant derivatives
of $A$ up to some specific order. We remark that in general Theorem \ref{compactness0}
fails to be true in the case $p=m$; in \cite{langer} on p.\ 227, Langer constructs a counterexample
in dimension $2$ by considering suitable inversions of a Clifford torus.
A similar result was shown by C.\ B.\ Ndiaye and R.\ Sch\"{a}tzle in \cite{ndiaye}, considering surfaces
with $L^2$-bounded second fundamental form that satisfy some additional hypotheses. Furthermore, the author showed in \cite{breuning1}
compactness of immersions with local Lipschitz representation.\\ \\
To prove Theorem \ref{compactness0}, we will first show a weak notion of convergence, the convergence in the sense of graph
systems. However, this does not directly imply the existence of a limit immersion
$f:M\rightarrow \R^n$. In \cite{langer}, in the case of surfaces, one
defines $M^j$ as limit manifold; here $j$ is a fixed large integer. Afterwards one constructs the
mappings $\phi^i:M^j\rightarrow M^i$ and shows, after passing to a subsequence, convergence to an immersion
$f:M^j\rightarrow \R^3$. Here we like to take a more systematic approach. We will
construct the limit manifold and immersion directly after having shown convergence of graph
systems. In order to do so, we shall take the limit graph system and define appropriate identifications;
this will enable us to recover the limit immersion by its image.
Only after that, we construct the diffeomorphisms $\phi^i$. This abstract construction
of the limit $f$ might be of its own interest for other applications. \\ \\
As a corollary of Theorem \ref{compactness0} we directly obtain:
\begin{cor}
Let $\mathfrak{F}$ be defined as in Theorem \ref{compactness0}.
Then there are only finitely many manifolds in $\mathfrak{F}$ up
to diffeomorphism.
\end{cor}
Next we prove a localized version
for smooth proper immersions admitting uniform $L^\infty$-bounds for the
second fundamental form $A$ and its covariant derivatives $\nabla^k A$. Here the manifolds on which the immersions
are defined are not required to be compact.
For a proper immersion $f:M\rightarrow \R^n$ with induced metric $g$ and volume
measure $\mu_g$ on $M$, let $\mu=f(\mu_g)$ be the Radon measure on $\R^n$ defined by
$\mu(E)=\mu_g(f^{-1}(E))$ for $E\subset\R^n$. Abbreviating we write $\|\cdot\|_{L^\infty(B_{R})}$ for the $L^\infty$-norm
on $f^{-1}(B_R)$, where $B_R\subset\R^n$ is the open ball of radius $R$ centered at the origin.
We obtain the following theorem: \\
\begin{theorem} \label{compactness1} (Compactness theorem for proper immersions) \\ \\
Let $f^i:M^i\rightarrow \R^n$ be a sequence of proper
immersions, where $M^i$ is an $m$-manifold without boundary and
$0 \in f^i(M^i)$. With
$\mu^i=f^i(\mu_{g^{i}}\!)$ assume
\beqn \mu^i(B_R) &\leq& C(R) \;\;\,\text{ for any } R>0, \label{comp1eq1} \\
\|\nabla^k A^i\|_{L^{\infty}(B_{R})} &\leq& C_k(R) \;\text{ for any } R>0 \text{ and }\,k \in
\N_0. \label{comp1eq2} \eeqn
Then there exists a proper immersion $f:M\rightarrow
\R^n$, where $M$ is again an $m$-manifold without boundary, such
that after passing to a subsequence there are
diffeomorphisms \beqnn \phi^i:U^i\rightarrow (f^i)^{-1}(B_i)\subset M^i, \eeqnn
where $U^i\subset M$ are open sets with $U^i\subset\!\subset U^{i+1}$
and $M=\bigcup_{i=1}^{\infty}U^i$, such that
$\|f^i\circ \phi^i-f\|_{C^{0}(U^i)}\rightarrow 0$,
and moreover $f^i\circ \phi^i\rightarrow f$ locally smoothly on $M$. \\ \\
Moreover, the immersion $f$ also satisfies (\ref{comp1eq1}) and (\ref{comp1eq2}), that is
$\mu(B_R) \leq C(R)$ and $\|\nabla^k A\|_{L^{\infty}(B_{R})}\linebreak \leq C_k(R)$.
\end{theorem} \vspace{5mm}
Again, the assumption $0\in f^i(M^i)$ can be weakened to $f^i(M^i)\cap K\neq \emptyset$ for
a fixed compact set $K\subset\R^n$. In contrast to the compact case, here the bound $\|A\|_{L^\infty(B_R)}\leq C(R)$
depends on the radius of the image. This explains the need of
some technical refinements which allow us to handle an increasing norm of the second fundamental form.
We like to remark that a similar result is shown by A.\ Cooper in \cite{cooper}, however the construction
of the diffeomorphisms $\phi^i$ is not carried out there
(see Remark \ref{cooprem} in this paper). Theorem \ref{compactness1}
has some important applications such as convergence proofs for geometric flows ---
for example for the mean curvature flow or the Willmore flow
(see e.g.\ \cite{baker1}, \cite{baker2}, \cite{huisken}, \cite{kuwert}, \cite{link}). \\ \\
As a corollary of Theorem \ref{compactness1} we prove
convergence of the corresponding measures: \\[-2mm]
\begin{cor} \label{measureconvergence} Let $f^i$ and $f$ be as in Theorem \ref{compactness1}
and let $\mu^i=f^i(\mu_{g^i})$, $\mu=f(\mu_g)$. Then
$\mu^i\rightarrow \mu$ in $C_c^0(\R^n)'$ as $i\rightarrow\infty$. \end{cor}
Finally, we will give some further generalizations of Theorem \ref{compactness1}.
In particular, in Corollary \ref{coropen}, we shall give a generalization
to proper immersions $f^i:M^i\rightarrow \Omega$ into an open subset $\Omega\subset \R^n$.
Along with this corollary, our theorems cover a wide range of situations one encounters in various applications.
\\ \\
\emph{Acknowledgement:} I would like to thank my advisor Ernst Kuwert for his support. Moreover
I would like to thank Manuel Breuning for proofreading
my dissertation \cite{breuningdiss}, where the results of this paper were established first.
\end{section}
\begin{section}{Local representation as a function graph} \label{locrepfg}
First, in Sections \ref{locrepfg} to \ref{convergenceimm}, we will show Theorem \ref{compactness0}.
After a rotation and a translation,
every immersion $f:M^m\rightarrow\R^{m+k}$ can locally be written as the graph
of a function $u:B_r\rightarrow \R^k$, where $B_r$ denotes an open ball in $\R^m$
of radius $r$. In this section, we like to work out the details
of such graph representations. First we have to introduce some notation: \\ \\
For $n=m+k$ let $G_{n,m}$ denote the Grassmannian
of (non-oriented) $m$-dimensional subspaces of $\R^n$. Unless stated otherwise, let
$B_{\varrho}$ always denote the open ball in $\R^m$
of radius $\varrho>0$ centered at the origin. \\ \\
Now let $M$ be an $m$-dimensional manifold without boundary and $f:M\rightarrow\R^n$
a $C^1$-immersion. Let $q\in M$ and let $T_qM$ be the tangent space at $q$. Identifying vectors
$X\in T_qM$ with $f_\ast X\in T_{f(q)}\R^n$, we may consider $T_qM$ as an $m$-dimensional
subspace of $\R^n$. Let $(T_qM)^\bot$ denote the orthogonal complement of
$T_qM$ in $\R^n$, that is \beqnn
\R^n=T_qM\oplus(T_qM)^\bot \eeqnn
and $(T_qM)^\bot$ is perpendicular to $T_qM$. In this manner
we may define a tangent and a normal map
\\ \parbox{12cm}{\beqnn \tau_f:M&\rightarrow& G_{n,m}, \hspace{3cm}\\
q&\mapsto& T_qM,  \\ \text{and} \hspace{6.3cm} \\
 \nu_f:M&\rightarrow& G_{n,k}, \\
q&\mapsto& (T_qM)^\bot. \eeqnn}   \hfill  \parbox{8mm}{\beqn \label{notiontangentspace} \eeqn
\vspace{0.2cm} \beqn \label{normalnotion} \eeqn} \\
Moreover, let $\pi_q^\top:\R^n\rightarrow T_qM$ and
$\pi_q^\bot:\R^n\rightarrow (T_qM)^\bot$ be the orthogonal projections onto
$T_qM$ and onto $(T_qM)^\bot$ respectively. \\ \\
First we like to consider immersions, that are already given as a graph. We like to begin with the
following trivial lemma:

\begin{lemma} \label{graphenabsch} Let $u,v:V\rightarrow \mathbb{R}^k$ be two
$C^1$-mappings, where $V \subset \mathbb{R}^m$ is open and convex
with $0 \in V$.
Moreover let $f,g:V \rightarrow \mathbb{R}^{m+k}$, $f(x)=(x,u(x))$, $g(x)=(x,v(x))$.
\begin{enumerate}[a)]
\item \label{graphtangente}The tangent space $\tau_f(x)$ is spanned by the vectors
$(e_1,\partial_1u(x)),\ldots, (e_m,\partial_mu(x))$.
\item \label{graphsteigung}If $u(0)=0$, then $|u(x)|\leq \|Du\|_{C^{0}(V)}|x|$.
\item \label{projektiontangente}Let $\zeta=(y,0) \in \mathbb{R}^m\times \mathbb{R}^k$. Then
$|\pi_x^\top(\zeta)|\geq
(1+\|Du(x)\|^2)^{-\frac{1}{2}}|y|$.
\item \label{projektionnormal}Let $\xi=(0,z) \in
\mathbb{R}^m\times \mathbb{R}^k$. Then
$|\pi_x^\bot(\xi)|\geq(1+\|Du(x)\|^2)^{-\frac{1}{2}}|z|$.
\end{enumerate}
\end{lemma}
The proof of the lemma is trivial and shall be omitted here. In the next lemma, we estimate the
$L^p$-norm of the second derivatives of $u$ from above by the supremum norm of the first derivative and the $L^p$-norm
of the second fundamental form:
\begin{lemma} \label{abschlp} For $B_r \subset \mathbb{R}^m$ and $n=m+k$, let $f \in
W^{2,p}(B_r,\mathbb{R}^n)$ be a mapping of the form $f(x)=(x,u(x))
\in \mathbb{R}^m \times \mathbb{R}^k$. If $\|Du\|_{C^{0}(B_{r})}<
\infty$, then we have the estimate \beg \|D^{2}u\|_{L^{p}(B_{r})}
\text{ }\leq \text{ } (1+\|Du\|_{C^{0}(B_{r})}^{2})^{\frac{3}{2}}\;
\|A(f)\|_{L^{p}(B_{r})}. \en
\end{lemma}
\textbf{Proof:} \\
Let $q \in B_r$. With Lemma \ref{graphenabsch}
\ref{projektionnormal}) we have \beqnn
|A_q(e_i,e_j)|&=&|\pi_q^\bot(\partial_{ij}f(q))|\\ &=&|\pi_q^\bot(0,\partial_{ij}u(q))|\\
&\geq& (1+\|Du(q)\|^2)^{-\frac{1}{2}}|\partial_{ij}u(q)| \\
&\geq&
(1+\|Du\|_{C^{0}(B_{r})}^{2})^{-\frac{1}{2}}|\partial_{ij}u(q)|.
\\\text{It follows} \hspace{3.5cm}\\ |\partial_{ij}u(q)|&\leq&
(1+\|Du\|_{C^{0}(B_{r})}^2)^{\frac{1}{2}}|A_q(e_i,e_j)| \\ &\leq&
(1+\|Du\|_{C^{0}(B_{r})}^2)^{\frac{1}{2}}|(e_i,\partial_iu(q))||(e_j,\partial_ju(q))|\|A(q)\|\hspace{1.5cm}
\\ &\leq& (1+\|Du\|_{C^{0}(B_{r})}^2)^{\frac{3}{2}}\|A(q)\|. \eeqnn
Integration yields the desired inequality. \hfill $\square$ \\ \\
\noindent
The following inequality is due to C.\ B.\ Morrey:
\begin{lemma} \label{absch1} Let $p>m$, $B_r\subset \mathbb{R}^m$ and $v \in (W^{1,p}\cap C^0)(B_r)$.
Then there is a universal constant $C=C(m,p)$, such
that for all $x \in B_r$ \beqn
\label{ungleichungmorrey}|v(x)-v(0)|\leq
Cr^{1-\frac{m}{p}}\|Dv\|_{L^{p}(B_{r})}. \eeqn
\end{lemma}
\textbf{Proof:} \\
For a proof see for instance \cite{alt}, p.\ 315, Theorem 8.11. The special case pointed out
on p.\ 317 in Remark 8.12 2) and 3)
is exactly (\ref{ungleichungmorrey}). \hfill
$\square$
\\ \\
With this lemma, we are able to estimate the supremum norm of the derivative from above
by the $L^p$-norm of the second derivatives:
\begin{lemma} \label{absch2} Let $p>m$, $B_r \subset \mathbb{R}^m$ and
$u \in (W^{2,p}\cap C^1)(B_r,\mathbb{R}^k)$. Let $u$ satisfy
$Du(0)=0$. Then there is a universal constant $C=C(m,k,p)$, such that \beqn
\|Du\|_{C^{0}(B_{r})} \leq
Cr^{1-\frac{m}{p}}\|D^{2}u\|_{L^{p}(B_{r})}. \eeqn
\end{lemma}
\textbf{Proof:}  \\
Using Lemma \ref{absch1}, with $Du \in (W^{1,p} \cap
C^0)(B_r,\mathbb{R}^{k \times m})$ the statement follows. \hfill $\square$ \\ \\
Next we like to explain, how an immersion can locally be written as a function graph.
The existence of such a graph representation is clear by the implicit function theorem.
However, for our purposes, we have to go more into detail.
First we need to introduce some more notation.
\\ \\
We call a mapping $A:\R^n\rightarrow \R^n$ a \emph{Euclidean isometry}, if
there is a rotation $R\in \mathbb{SO}(n)$ and a translation $T\in\R^n$, such that
$A(x)=Rx+T$ for all $x\in \R^n$. \\ \\
\label{defeuclisom}For a given point $q \in M$ let $A_q: \mathbb{R}^n\rightarrow
\mathbb{R}^n$ be a Euclidean isometry,
which maps the origin to
$f(q)$, and the subspace $\mathbb{R}^m\times\{0\}\subset
\mathbb{R}^m \times \mathbb{R}^k$ onto $f(q)+\tau_f(q)$.
Let $\pi:\R^n\rightarrow \R^m$ be the standard projection onto the
first $m$ coordinates. \\ \\
Finally let $U_{r,q}\subset M$ be the $q$-component
of the set \label{qcomponentof} $(\pi \circ A_q^{-1} \circ f)^{-1}(B_r)$.
Although the isometry $A_q$ is not uniquely determined, the set
$U_{r,q}$ does not depend on the choice of $A_q$. \\ \\
We come to the central definition (as first defined in \cite{langer}):
\begin{defi}
An immersion $f:M\rightarrow \mathbb{R}^n$ is called an $(r,\alpha)$-immersion,
if for each point $q\in M$ the set $A_q^{-1}\circ f(U_{r,q})$
is the graph of a differentiable function $u:B_r\rightarrow \mathbb{R}^k$ with
$\|Du\|_{C^{0}(B_{r})}\leq \alpha$. \end{defi} \vspace{3mm}
Here, for any $x\in B_r$ we have $Du(x)\in \R^{k\times m}$. In order to define the
$C^0$-norm for $Du$, we have to fix a matrix norm for $Du(x)$. Let us agree upon
\vspace{-1.5mm}
\beqnn \|A\|\:=\:\Biggl(\hspace{0.5mm}\sum_{j=1}^m|a_j|^2\Biggr)^{\!\frac{1}{2}} \eeqnn
for $A=(a_1,\ldots,a_m)\in \R^{k\times m}$. For this norm we have
$\|A\|_{\text{op}}\:\leq\: \|A\|$ for any $A\in \R^{k\times m}$
and the operator norm $\|\cdot\|_{\text{op}}$.
Hence the bound $\|Du\|_{C^0(B_r)}\leq \alpha$ directly implies that $u$
is $\alpha$-Lipschitz (and all estimates in the previous lemmas are true).
Moreover the norm $\|Du\|_{C^{0}(B_{r})}$ does not depend on the choice of the
isometry $A_q$. \\ \\
For given $\alpha>0$, we would like to give an estimate for the admissible size of the radius $r$,
such that a given immersion is an $(r,\alpha)$-immersion.
Here the admissible size of $r$ only depends
on the $L^p$-norm of the second fundamental form:
\begin{theorem} Let $p>m$ and $0<\alpha \leq1$. Then there exists a universal constant
$c=c(m,k,p)>0$ such that every immersion $f\in W^{2,p}(M,\R^n)$ on a compact $m$-manifold
$M$ is an $(r,\alpha)$-immersion for all $r>0$ with \vspace{-1.2mm} \beqn \label{sizegraphradius}
r^{1-\frac{m}{p}}\:\leq\: c\alpha\|A(f)\|_{L^{p}(M)}^{-1}. \eeqn
\end{theorem}
\textbf{Proof:}
The proof of the $2$-dimensional case in \cite{langer} also applies to the higher dimensional case,
where we use Lemmas \ref{abschlp} and \ref{absch2}. \hfill $\square$ \\ \\
In the previous lemma, $k$ denotes the codimension. The following lemma
(Lemma 3.1 in \cite{langer}) is crucial for the proof
and will also be needed (in a variation) for the noncompact case. The proof of Langer
carries over to our situation: \\[-2.5mm]
\begin{lemma} \label{inclusiondipl}
Let $f:M\rightarrow \R^n$ be an $(r,\alpha)$-immersion and $p,q \in M$.
\begin{itemize}
\item[a)] If $p \in U_{r,q}$, then
$|f(p)-f(q)|\leq (1+\alpha^2)r$.
\item[b)] If $\alpha^2<\frac{1}{3}$ and
$U_{\frac{r}{4},q}\cap U_{\frac{r}{4},p}\neq \emptyset$, then
$U_{\frac{r}{4},p}\subset U_{r,q}$.
\end{itemize}
\end{lemma} \vspace{3mm}
For $(r,\alpha)$-immersions $f:M\rightarrow \R^n$ we define the notion of a $\delta$-net:
\begin{defi} Let $Q=\{q_1,\ldots,q_s\}$ be a finite set of points in $M$ and let
$0<\delta<r$. We say that $Q$ is a $\delta$-net for $f$ if $M=\bigcup_{j=1}^s U_{\delta,q_{j}}$.
\end{defi}
The number of elements of a $\delta$-net can be bounded from above:
\begin{lemma} \label{boundindiss} Assume $\alpha^2<\frac{1}{3}$ and $0<\delta<r$.
Then every $(r,\alpha)$-immersion
$f:M\rightarrow \R^n$ admits a $\delta$-net with at most
$\left(\frac{4}{\delta}\right)^m \vol(M)$ points. \end{lemma}
\textbf{Proof:}\\
The proof is the same as in the $2$-dimensional case, see Lemma 3.2 in \cite{langer}. Note that one could even
derive the bound $\left(\frac{4}{\delta}\right)^m(\mathcal{L}^m(B_1))^{-1}\,\text{vol}(M)$. \hfill $\square$
\end{section}
\begin{section}{Convergence of graph systems}
In the previous section we have seen, how any immersion in $\mathfrak{F}$ can be written locally on sets $U_{r,q}$
as the graph of a function. The notion of a $\delta$-net yields a cover
of each manifold with such kind of sets. This is the starting point for the notion of graph
systems, and for convergence of such systems. \\ \\
First we like to explain how to represent an immersion in $\mathfrak{F}$ as a system of graphs.
For that we define the space of graph systems with $s$ elements by \\
\parbox{12cm}{\beqnn \mathfrak{G}^s=\{(A_j,u_j)_{j=1}^s: &A_j&:\:
\mathbb{R}^n\rightarrow \mathbb{R}^n \text{ is a Euclidean
isometry,} \\ &u_j& \in
W^{2,p}(B_r,\mathbb{R}^k)\}. \eeqnn }  \hfill\parbox{8mm}{\beqn \label{grundraum} \eeqn} \\
Every Euclidean isometry $A:\R^n\rightarrow \R^n$ splits uniquely into a rotation
$R\in \mathbb{SO}(n)$ and a translation
$T\in \R^n$, such that $A(x)=Rx+T$ for all $x\in \R^n$.
If $\|\cdot\|_{\text{op}}$ denotes the operator norm and if
$\Gamma=(A_j,u_j)_{j=1}^s\in \mathfrak{G}^s$,
$\tilde{\Gamma}=(\tilde{A}_j,\tilde{u}_j)_{j=1}^s\in \mathfrak{G}^s$, we set
\\ \parbox{12cm}{ \beqnn \mathfrak{d}(\cdot,\cdot): \mathfrak{G}^s\times
\mathfrak{G}^s&\rightarrow& \mathbb{R}, \\
\mathfrak{d}(\Gamma,\tilde{\Gamma})&=&\sum
\limits_{j=1}^{s}(\|R_j-\tilde{R}_j\|_{\text{op}}+|T_j-\tilde{T}_j|+\|u_j-\tilde{u}_j\|_{C^{1}(B_{r})}).
\eeqnn} \hfill \parbox{8mm}{\beqn \label{metrikgraph1} \eeqn} \\
This makes $(\mathfrak{G}^s,\mathfrak{d})$ a metric space. \\ \\
Now let $f:M\rightarrow \R^n$ be an $(r,\alpha)$-immersion and $Q=\{q_1,\ldots,q_s\}$ a
$\delta$-net for $f$ with $s$ elements. To each $q_j\in Q$ we may assign a neighborhood
$U_{r,q_{j}}$, a Euclidean isometry $A_j$, and a $C^1$-function $u_j:B_r\rightarrow \R^k$
as described above. Hence, to given $f$, $r$ and $Q$, we may assign a graph system \beqnn
\Gamma=\Gamma(f)=(A_j,u_j)_{j=1}^s \in \mathfrak{G}^s. \eeqnn
The isometries $A_j$ and functions $u_j$ are not uniquely determined, but we always have
$u_j(0)=0$ and $Du_j(0)=0$. \\ \\
For any $j\in\{1,\ldots,s\}$ we finally set
$Z(j):=\{1\leq k \leq s:U_{\delta,q_{j}}\cap U_{\delta,q_{k}}\neq \emptyset\}$. \\ \\ \\[-4.4mm]
With the preceding notations we are able to define a notion of convergence for
graph systems: \\[-2.5mm]
\begin{defi}[Convergence in the sense of graph systems] $ $ \\
Let a sequence $f^i: M^i\rightarrow \mathbb{R}^n$ of immersions be given.
We say $f^i$ is convergent in the sense of graph systems,
if there are fixed \beqn \alpha,r,\delta > 0 \text{ with }
r > \delta, \text{ and } s \in \mathbb{N}, \nonumber \eeqn
such that the following properties are satisfied: \begin{itemize}
\item[$-$] Each $f^i$ is an $(r,\alpha)$-immersion.
\item[$-$] For each $f^i$ there exists a $\delta$-net with $s$
points, for which the following holds: \begin{itemize}
\item[$\bullet$] $Z^i(j)=Z(j)$ for fixed sets $Z(j)$ independent of \,$i$.
\item[$\bullet$] There exists a system $\Gamma \in \mathfrak{G}^s$, such that
the graph systems $\Gamma^i$ corresponding to $f^i$ converge in
$(\mathfrak{G}^s,\mathfrak{d})$ to $\Gamma$.
\end{itemize} \end{itemize}
\end{defi}
\vspace{2.5mm}
The following statement is true:
\begin{theorem} \label{graphsystemconv} Every sequence in $\mathfrak{F}$ admits a subsequence
that converges in the sense of graph systems. \end{theorem}
\textbf{Proof:} \\
Using the results above, the proof of Theorem 3.3 on p.\ 228 in \cite{langer} carries over
to the higher dimensional case. \hfill $\square$ \\ \\
Here, we only require a graph system $\Gamma$ as limit, but not an immersion $f$. Actually we could say,
that any sequence in $\mathfrak{F}$ admits a subsequence that is Cauchy in the sense of graph systems.
In the next section we will show completeness in the sense that there exists an immersion
$f$ with $\Gamma=\Gamma(f)$.
\end{section}
\begin{section}{Construction of the limit manifold and immersion} \label{constructionsection}
In Theorem \ref{graphsystemconv}, for a given sequence of immersions $f^i$ in $\mathfrak{F}$ we have found
a subsequence, that converges in the sense of graph systems to a limit system $\Gamma$. However, it is not clear
whether $\Gamma$ is the graph system of an immersion $f:M\rightarrow \R^n$ on a compact manifold $M$. In this
section we like to show, that this is the case. \\ \\
First we would like to construct the limit manifold $M$. We start with a sequence of $(r,\alpha)$-immersions, convergent
in the sense of graph systems, with $\alpha^2<\frac{1}{3}$, $\delta=\frac{r}{16}$,
$\frac{\delta}{10}$-nets $Q^i=\{q_1^i,\ldots,q_s^i\}$ with $s$ elements, intersection
sets $Z(j)=\{1\leq
k \leq s: U_{\delta,q_j^i}^i\cap U_{\delta,q_k^i}^i\neq \emptyset\}$ which are
independent of $i$, limit isometries $A_j$ and limit functions
$u_j:B_r\rightarrow \R^k$. Here we have to use $\delta/10$-nets
and not only $\delta$-nets; this is in particular needed in the proof of Lemma \ref{eigenschaft1}.
To simplify the notation, for $0<\varrho\leq r$ we set
$U_{\varrho,j}^i\,:=\,U_{\varrho,q_j^i}^i$.
\\ \\
For the open ball $B_\delta\subset \R^m$ we set $B_\delta^j=B_\delta\times\{j\}$.
This makes $\bigcup_{j=1}^sB_\delta^j$ a disjoint union. We endow $\bigcup_{j=1}^sB_\delta^j$
with the \emph{topology of the disjoint union,} which is defined as follows: A subset
$U\subset \bigcup_{j=1}^sB_\delta^j$ is open if and only if $U\cap B_\delta^j\subset B_\delta$
is open for every $j$.
\\ \\ We define a relation $\sim$ on $\bigcup_{j=1}^{s}B_\delta^j$.
For $(x,j), (y,k) \in \bigcup_{l=1}^{s}B_\delta^l$
we set \beqn (x,j)\sim (y,k) \;
\Leftrightarrow \; [k \in Z(j) \text{ and }
A_j(x,u_j(x))=A_k(y,u_k(y))]. \eeqn
\begin{lemma} \label{equivalence}
The relation $\sim$ is an equivalence relation. \end{lemma}
\textbf{Proof:} \\
Obviously the relation $\sim$ is reflexive and symmetric.
Now let $(x,j)\sim(y,k)$, $\;
(y,k)\sim(z,l)$ for $(x,j),(y,k),(z,l) \in
\bigcup_{\nu=1}^{s}B_\delta^\nu$.
As $k \in Z(j)$, $l \in Z(k)$, we have $U_{\delta,j}^{i}\cap
U_{\delta,k}^{i}\neq \emptyset$, $U_{\delta,k}^{i}\cap
U_{\delta,l}^{i} \neq \emptyset$. Using
Lemma \ref{inclusiondipl} b) twice yields $U_{\delta,l}^{i}\;
\subset \; U_{r,j}^{i}$. Moreover there is exactly one
$\xi^{i} \in U_{\delta,l}^{i}$ with $f^{i}(\xi ^{i})=A_{l}^{i}(z,u_{l}^{i}(z))$.
By the definition of $\sim$ it follows $A_j(x,u_j(x))=A_l(z,u_l(z))$,
which means together with the graph convergence
$f^i(\xi ^i) \rightarrow
A_j(x,u_j(x)) \text{ as } i \rightarrow \infty$. We
define $x^i := \pi\circ (A_{j}^{i})^{-1} \circ f^{i}(\xi^i)$.
As $U_{\delta,l}^{i}\subset U_{r,j}^{i}$, we have $\xi^i\in
U_{r,j}^i$ and hence $x^i \in B_r$ and
$f^i(\xi^i)=A_j^i(x^i,u_j^i(x^i))$. With the convergence of
$f^i(\xi^i)$ it follows $A_j^i(x^i,u_{j}^{i}(x^i)) \rightarrow
A_j(x,u_j(x))$ as $i\rightarrow \infty$. As $A_j^i
\rightarrow A_j$ for $i\rightarrow \infty$ (in the sense \nolinebreak of
(\ref{metrikgraph1})) it follows $\pi\circ(A_j^i)^{-1}\circ
A_j(x,u_j(x))\rightarrow x$, and, by using the triangular inequality,
$x^i \rightarrow x$ in \nolinebreak the ball $B_r$.
In particular we have $x^i \in B_{\delta}$ for $i$ sufficiently large, and hence $\xi^i\in U_{\delta,j}^i$.
We deduce $U_{\delta,j}^i\cap U_{\delta,l}^i \neq \emptyset$ (which is then automatically satisfied
for all $i$) and hence $l \in
Z(j)$. This shows transitivity. \hfill
$\square$ \\ \\
This enables us to define the limit manifold.
As set, $M$ is defined to be the quotient space
\beqn M=\Biggl(\bigcup_{j=1}^{s}B_\delta^j\Biggr)/\sim \;\:, \eeqn
resulting from the equivalence relation
of above. Let $M$ be endowed with the quotient topology. \label{constructionM}
For $(x,k) \in
\bigcup_{j=1}^{s}B_\delta^j$, let $[(x,k)]$ denote the
corresponding equivalence class. Let $P$ denote the canonical projection
from $\bigcup_{j=1}^{s}B_\delta^j$ onto $M$,
and $P_j$ the restriction
$P|B_\delta^j:\!B_\delta^j\!\rightarrow\! P(B_\delta^j)$. We can consider
$P_j$ as a mapping defined on $B_\delta$. We note that
$P$ is injective on $B_\delta^j$, and in particular $P_j$ invertible.
For $V\subset B_\delta$ we set $V^j=V\times\{j\}$.
For $V \subset B_{\delta}$ open we define \\[-6mm] \\ \parbox{12cm}{\beqnn \varphi_{_V}^{j}: P(V^j)
&\rightarrow& V,
\\ \text{[}(x,j) \text{]}&\mapsto& P_{j}^{-1}([(x,j)]) \in V, \eeqnn} \hfill
\parbox{8mm}{\beqn \eeqn} \\
which yields a well-defined mapping. Finally we denote the
set of all such mappings by $\mathfrak{A}$, that is\,
$\mathfrak{A}=\{\varphi_{_W}^{k}:\;1\leq k \leq s,\: W\subset B_\delta\, \text{ open}\}$.
To simplify the notation, we will often identify sets $V^j$ with
$V$, and elements $(x,j)$ with $x$ (as already done above).
\begin{lemma} \label{quotientopen}
The quotient projection $P$ is open. \end{lemma}
\textbf{Proof:} \\
Let $V \subset \bigcup_{j=1}^{s}B_\delta^j$ be open. We have to show,
that $P^{-1}(P(V))$ is open. For that let $x \in P^{-1}(P(V))$.
Then $x \in B_{\delta}^{j}$ for a $j \in
\{1,\ldots,s\}$. We show the existence of an open neighborhood $U
\subset B_{\delta}^j$ of $x$ with $U \subset P^{-1}(P(V))$. \\ \\
It holds $x\sim y$ for a $y \in V$ and moreover $y \in B_{\delta}^k$
for a $k \in Z(j)$. Now consider $\psi:
B_\delta \rightarrow \mathbb{R}^m$, $z\mapsto \pi\circ A_k^{-1}\circ
A_j(z,u_j(z))$. As $x\sim y$, we have $\psi(x)=y$. As $V$ is open, there
is an open neighborhood $W \subset V$ of $y$ with $W\subset
B_{\delta}^k$. As $\psi$ is continuous, $\psi^{-1}(W)$
is an open neighborhood of $x$. \\ \\
We show that every point $z \in \psi^{-1}(W)$ is equivalent to a point
in $W$, which implies the statement. For
every $i$ there is exactly one $\xi^i \in U_{\delta,j}^i$ with
$f^i(\xi^i)=A_j^i(z,u_j^i(z))$. As $k \in Z(j)$, with Lemma
\ref{inclusiondipl} b) it holds $U_{\delta,j}^i\subset U_{4\delta,k}^i$.
Hence for every $i$ there is a $w^i \in B_{4\delta}$ with
$A_k^i(w^i,u_k^i(w^i))=A_j^i(z,u_j^i(z))$. For $i\rightarrow
\infty$ we have $A_j^i(z,u_j^i(z))\rightarrow A_j(z,u_j(z))$, and
$A_k^i\rightarrow A_k$, $u_k^i\rightarrow u_k$ (in the sense of
(\ref{metrikgraph1})) and for a subsequence $w^i\rightarrow w$
for a $w \in \overline{B_{4\delta}}$. Using the triangular inequality, we deduce
$A_k(w,u_k(w))=A_j(z,u_j(z))$.
Hence $w=\psi(z) \in W\subset B_\delta$ and $(z,j)\sim
(w,k)$, which proves the lemma. \hfill $\square$ \\[-2.25mm]
\begin{lemma} \label{MisHausdorff}
The space $M$ is a second countable Hausdorff space. \end{lemma}
\textbf{Proof:}\\
We first show that $M$ is Hausdorff.
Let $p,q \in M$ with $p \neq q$. Then there are $j,k \in
\{1,\ldots,s\}$ with $p \in P(B_\delta^j)$, $q \in P(B_\delta^k)$.
If $k \notin Z(j)$, then $P(B_\delta^j)$ and $P(B_\delta^k)$
are disjoint open neighborhoods. \\ \\
Now let us assume $k \in Z(j)$. Then there are $x \in B_\delta^j$, $y
\in B_\delta^k$ with $p=P(x)$, $q=P(y)$. It follows $A_j(x,u_j(x)) \neq
A_k(y,u_k(y))$, as otherwise $p=q$. We define a mapping
$\gamma: B_\delta\times B_\delta\rightarrow \mathbb{R}$,
$(v,w)\mapsto |A_j(v,u_j(v))-A_k(w,u_k(w))|$. Hence
$\gamma(x,y) > 0$. As $\gamma$ is continuous, there are open neighborhoods
$V,W$ of $x,y$ with
$\gamma(V\times W)\subset (0,\infty)$. Using that the projection $P$ is open,
$P(V^j)$ and $P(W^k)$ are disjoint open neighborhoods of $p$ and $q$. \\ \\
Next we like to show that $M$ is second countable. Let $\mathfrak{B}$ be a countable
basis of $\bigcup_{j=1}^{s}B_\delta^j$. As the projection $P$ is open, $\{P(B): B \in \mathfrak{B}\}$
is a countable basis of $M$. \hfill $\square$ \\
\begin{lemma} \label{diffatlas}
The set $\mathfrak{A}$ is a differentiable atlas on $M$. \end{lemma}
\textbf{Proof:} \\
First we note that $M$ is covered by the sets $P(V^j)$.
Furthermore, every
$\varphi_{_V}^j:P(V^j)\rightarrow V$ is a bijective mapping between
open sets with inverse mapping $P_j$
(more precisely $(\varphi_{_V}^j)^{-1}: V\rightarrow P(V^j)$ with
$(\varphi_{_V}^j)^{-1}(x)=P_j(x)$). The
quotient projection $P$ is open by Lemma \ref{quotientopen},
and certainly continuous. But continuous, open, bijective mappings are homeomorphisms.
Hence $M$ is locally Euclidean.
\\ \\ It remains to show differentiability of the coordinate changes.
For the charts
$\varphi_{_V}^j, \varphi_{_W}^k$, the coordinate change is given by \beqnn
\varphi_{_V}^j\circ
(\varphi_{_W}^k)^{-1}: \;\;\;\;\varphi_{_W}^k(P(V^j)\cap P(W^k))&\rightarrow& \varphi_{_V}^j(P(V^j)\cap P(W^k)), \\
x\;\;\;\;\;\;\;\;\;\;\;\;\;&\mapsto& \pi\circ A_{j}^{-1}\circ
A_{k}(x,u_k(x)). \eeqnn
But this is a composition of smooth mappings; hence $\varphi_{_V}^j\circ(\varphi_{_W}^k)^{-1}$ is smooth.
\hfill $\square$ \\ \\
Let us summarize our results:
\begin{theorem} The topological space $M$ is Hausdorff with countable basis
and $\mathfrak{A}$ is a differentiable atlas on
$M$. Hence $(M,\mathfrak{A})$ induces uniquely the structure of a differentiable manifold.
\end{theorem}
Finally we show compactness of $M$:
\begin{lemma} \label{compact} The limit manifold $M$ is
compact.
\end{lemma}
\textbf{Proof:} \\
For the proof we already use Lemma \ref{eigenschaft1}. By this
we have $M=\bigcup_{j=1}^{s}P(\overline{B_{\delta/2}^{j}})$. As the quotient projection
is continuous, with the compactness of
$\overline{B_{\delta/2}^{j}}$ the statement follows.
\hfill $\square$ \\ \\ \\
Now we define the limit immersion:
\\
\parbox{12cm}{
\beqnn f:\;\;\;\;\; M\;\;\;&\rightarrow& \mathbb{R}^n, \\
\text{[}(x,j)\text{]}
&\mapsto& A_j(x,u_j(x)). \eeqnn} \hfill \parbox{8mm}{\beqn \eeqn} \\
If $(x,j)\sim(y,k)$, by the definition of $\sim$
we have $A_j(x,u_j(x))=A_k(y,u_k(y))$. Hence $f$ is
well-defined. Moreover $f$ admits the local representation $x\mapsto
A_j(x,u_j(x))$ for $x \in B_\delta$, which implies that $f$ is an
immersion. Finally we note that the limit system $(A_j,u_j)_{j=1}^s$
of the graph convergence is the graph system of an immersion. \\ \\ \\
The following lemmas are associated with the construction of the limit manifold above.
All statements are needed only for technical reasons and will be required
for the construction of the mappings $\phi^i$ in the next section, and in particular
for showing injectivity of these mappings.
Additionally, Lemma \ref{eigenschaft1} is required in the proof of Lemma \ref{compact},
stating that $M$ is compact. \\ \\
By the definition of $M$, we have
$M=\bigcup_{j=1}^{s}P(B_\delta^j)$. The following lemma says,
that there even exists a much finer cover:
\begin{lemma} \label{eigenschaft1}
It holds $M=\bigcup_{j=1}^{s}P(B_{\frac{\delta}{6}}^j)$.
\end{lemma}
\textbf{Proof:} \\
Let $q \in M$ be an arbitrary point. Then there is a $j\in \{1,\ldots,s\}$ and
an $x \in B_\delta^j$ with $P(x)=q$. It follows
$f(q)=A_j(x,u_j(x))=\lim_{i\rightarrow \infty} A_j^i(x,u_j^i(x))$.
Moreover, there are $\xi^i\in U_{\delta,j}^{i}$ with
$f^i(\xi^i)=A_j^i(x,u_j^i(x))$. As the sets
$Q^i=\{q_1^i,\ldots,q_s^i\}$ are $\frac{\delta}{10}$-nets
for $f^i$, there are $j^i\in \{1,\ldots,s\}$ with $\xi^i \in
U_{\delta/10,j^{i}}^{i}$. After passing to a subsequence, we may assume
$j^i=k$ independent of $i$.
Then there are $y^i \in B_{\delta/10}$ with
$f^i(\xi^i)=A_{k}^i(y^i,u_{k}^i(y^i))$. A subsequence of $y^i$
converges to $y \in \overline{B_{\delta/10}}\subset
B_{\delta/6}$. As $f^i(\xi^i)\rightarrow f(q)$,
$A_{k}^i\rightarrow A_{k}$ and $u_{k}^i\rightarrow u_{k}$ for
$i\rightarrow \infty$, we have $f(q)=A_{k}(y,u_{k}(y))$. As $\xi^i\in
U_{\delta,j}^i$, $\xi^i \in U_{\delta/10,k}^i$, we have $U_{\delta,j}^i
\cap U_{\delta/10,k}^i\neq \emptyset$, and all the more
$U_{\delta,j}^i\cap U_{\delta,k}^i \neq \emptyset$. This implies $k \in
Z(j)$, and moreover $(x,j)\sim (y,k)$. It follows $q \in
P(B_{\delta/6}^{k})$. \hfill $\square$
\\ \\
The next statement is the analogue to Lemma \ref{inclusiondipl}
 b) for the limit immersion:
\begin{lemma} \label{eigenschaft2}
If $P(B_{\frac{\delta}{4}}^j)\cap P(B_{\frac{\delta}{4}}^k) \neq
\emptyset$, then $P(B_{\frac{\delta}{4}}^k)\subset
P(B_{\delta}^j)$. \end{lemma}
\textbf{Proof:} \\
The proof of Lemma \ref{inclusiondipl} carries over to the limit immersion. \hfill $\square$
\\ \\ Analogous to the sets $Z(j)$, we define intersection sets for a finer cover of
$M^i$ by \beqnn
\tilde{Z}^i(j)=\{1\leq k \leq s: U_{\frac{\delta}{5},j}^i \cap
U_{\frac{\delta}{5},k}^i \neq \emptyset \}. \eeqnn
Passing to a subsequence, again we may assume
$\tilde{Z}^i(j) = \tilde{Z}(j)$ independent of $i$. \\ \\
The relation $P(B_\delta^j)\cap P(B_\delta^k) \neq \emptyset$ implies $k
\in Z(j)$; however, in general $P(B_\delta^j)\cap
P(B_\delta^k)
= \emptyset$ does not imply $k \notin Z(j)$. Instead the following
statement holds (where the numbers are adapted to the situation in the
next section):
\begin{lemma} \label{eigenschaft3}
If $P(B_{\frac{\delta}{4}}^j)\cap P(B_{\frac{\delta}{4}}^k) =
\emptyset$, then $k \notin \tilde{Z}(j)$. \end{lemma}
\textbf{Proof:} \\
Let $P(B_{\delta/4}^j)\cap P(B_{\delta/4}^k)=\emptyset$.
Suppose $k \in \tilde{Z}(j)$. Then for every
$i$ there is a $\xi^i \in U_{\delta/5,j}^i\cap U_{\delta/5,k}^i$.
Moreover, the points $f^i(\xi^i)$
lie in a ball of fixed radius. Hence there is a subsequence and an
 $x \in \mathbb{R}^n$ with
$f^i(\xi^i)\rightarrow x$ as $i \rightarrow \infty$. With the graph
convergence and by arguments as in Lemma
\ref{eigenschaft1}, we have $x=A_j(y,u_j(y))=A_k(z,u_k(z))$ with
$y,z \in \overline{B_{\delta/5}}\subset B_{\delta/4}$. As $k \in
\tilde{Z}(j)$, we surely have $k \in Z(j)$. It follows
$P(B_{\delta/4}^j)\cap P(B_{\delta/4}^k) \neq \emptyset$, contrary to our assumption.
\hfill $\square$
\end{section}
\begin{section}{Reparametrization of the immersions} \label{secrepimm}
We like to construct the reparametrizations $\phi^i:M\rightarrow M^i$.
This is done by a kind of projection from the limit surface onto each of the
surfaces $f^i$. \\ \\
Our starting point is a sequence of $(r,\alpha)$-immersions
$f^i:M^i\rightarrow \R^n$ in $\mathfrak{F}$, which converges in the sense
of graph systems to a limit immersion $f:M\rightarrow \R^n$. Here
we require \label{alphasquare} $\alpha^2 \leq \frac{1}{10}$.
We will define the projection locally, using charts
$\varphi_j:P(B_\delta^j)\rightarrow B_\delta$.
By such a chart, we shall often tacitly identify the set $P(B_\delta^j)$ with
the ball $B_\delta$. \\ \\
Let $A_j$ and $A_j^i$ denote the isometries of the previous sections
corresponding to $f$ and $f^i$ respectively.
As the following constructions are invariant under translations and rotations,
we may assume $A_j=\text{Id}_{\R^n}$ and replace $A_j^i$ by $A_j^{-1}\circ A_j^i$. \\ \\
Then $f(P(B_\delta^j))$ is the graph of a function $u_j:B_\delta\rightarrow \R^k$
with $u_j(0)=0$, $Du_j(0)=0$. The set $f^i(U_{r,j}^i)$ is the graph of a function
$u_j^i:B_r\rightarrow \R^k$, however translated and rotated relatively to the limit immersion
by $A_j^{-1}\circ A_j^i$. But actually $A_j^i\rightarrow A_j$
as $i\rightarrow\infty$ in the sense of the metric (\ref{metrikgraph1}).
Hence the translation and rotation $A_j^{-1}\circ A_j^i$ gets arbitrarily small
relative to $f$ as $i\rightarrow\infty$. Hence we may assume that also
$f^i(U_{r,j}^i)$ is the graph of a function on a subset of
$\R^m\subset \R^m\times \R^k$, which shall be denoted in the following by
$\tilde{u}_j^i$. \\ \\
Furthermore for all $\varrho$ with $0<\varrho<r$ there is an $N\in \N$,
such that for all $i>N$ \beqn \label{positionimmersion1}
\{(x,\tilde{u}_j^i(x)):x\in B_{r-\varrho}\}\subset f^i(U_{r,j}^i). \eeqn \\[-5mm]
This is the situation represented in Figure \ref{posimmersions}. \\[-5mm]
\setlength{\unitlength}{1cm}
\begin{picture}(10,6) \put(1.5,0){\put(-0.6,0){\line(1,0){11.2}} \put(1.2,-0.1){(}
  \put(8.8,-0.1){)} \put(-0.65,-0.1){(} \put(10.5,-0.1){)}
  \put(-0.1,-0.1){(} \put(10,-0.1){)} \put(9.5,-0.6){$B_{r-\varrho}$}
  \put(5,-0.4){0} \put(8.7,-0.6){$B_\delta$} \put(10.4,-0.6){$B_r$}
  \put(11,1.9){$f$ (limit immersion)}

  \qbezier[20](-0.6,1.2)(0.3,0.8)(1.2,0.6)
  \qbezier(1.2,0.6)(3.5,0)(5,0) \qbezier(5,0)(7.4,0)(8.2,0.5)
  \qbezier(8.2,0.5)(8.8,0.9)(8.8,0.9)
  \qbezier[20](8.8,0.9)(9.9,1.7)(10.6,2.1)

  \qbezier[200](-1.1,3)(5,1.8)(10.3,4)
  \qbezier[200](-0.4,2.4)(6,0.8)(10.8,3.4)
  \qbezier[200](-0.8,2.3)(6.7,-0.4)(10.6,2.7) \put(11,4.05){$f^1(U_{r,j}^1)$}
  \put(11,3.45){$f^2(U_{r,j}^2)$} \put(11,2.85){$f^3(U_{r,j}^3)$} \put(11.1,2.3){$\vdots$}

\qbezier[80](-0.045,0)(-0.045,2.1)(-0.045,4.2)
\qbezier[80](10.095,0)(10.095,2.1)(10.095,4.2)}
\end{picture} \\ \\
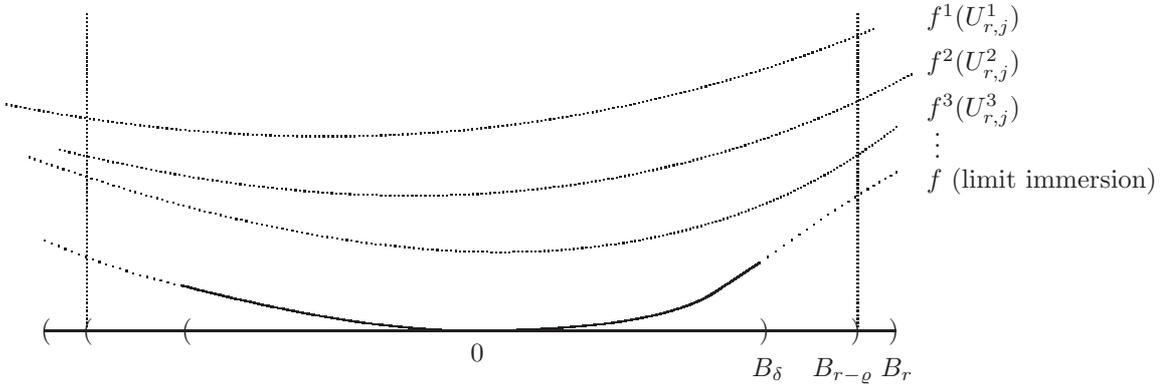
\begin{fig} \label{posimmersions} Position of the immersions $f^i$ relative to the limit immersion.
Note that the figure is not true to scale, because in the proof we have $r=16\delta$. \end{fig}
As $\|Du_j^i\|_{C^0(B_r)}\leq \alpha$, we surely may assume
$\|D\tilde{u}_j^i\|_{C^0(B_{r-\varrho})}\leq 2\alpha$ for $i$
sufficiently large.
Moreover, by the graph convergence, for any $\varepsilon>0$ we have
$|\tilde{u}_j^i(0)|<\varepsilon$ for $i$ large.
\\ \\
Finally we like to simplify notation. All the following considerations are
performed locally on $P(B_\delta^j)$. We will fix the index $j$ and suppress it in the notation.
Hence we shall write for example $u$ instead of $u_j$, and $\tilde{u}^i$
instead of $\tilde{u}_j^i$. \\ \\
If the limit immersion is sufficiently smooth, it is possible to project into the normal direction.
However, if $f$ is not $C^2$, in general this is not possible. Without
an $L^\infty$-bound for the second fundamental form we might have a local concentration of
curvature. In this case, projecting into the normal direction will not lead to injective mappings $\phi^i$
(see Figure \ref{normalconc}).

\begin{picture}(10,4.5) \put(0.5,0){ \put(10.8,2){$f$}

  \qbezier[20](-0.6,1.2)(0.3,0.8)(1.2,0.6)
  \qbezier(1.2,0.6)(3.5,0)(5,0) \qbezier(5,0)(7.4,0)(8.2,0.5)
  \qbezier(8.2,0.5)(8.8,0.9)(8.8,0.9)
  \qbezier[20](8.8,0.9)(9.9,1.7)(10.6,2.1)

  \put(7.25,0.17){\line(-1,6){0.1644}}
  \put(7.5,0.23){\line(-1,5){0.19612}}
  \put(7.7,0.3){\line(-1,4){0.24254}}
  \put(8.1,0.45){\line(-5,6){0.64018}}
  \put(7.9,0.35){\line(-1,3){0.31623}}
  \put(7.985,0.4){\line(-2,5){0.37139}}
  \put(8.07,0.42){\line(-3,5){0.5145}}
  \put(8.185,0.5){\line(-4,5){0.6247}}
  \put(8.27,0.55){\line(-3,4){0.6}}
  \put(8.47,0.7){\line(-2,3){0.5547}}
  \put(8.67,0.83){\line(-3,5){0.5145}}

  \put(7.25,0.17){\line(1,-6){0.1644}}
  \put(7.5,0.23){\line(1,-5){0.19612}}
  \put(7.7,0.3){\line(1,-4){0.24254}}
  \put(8.1,0.45){\line(5,-6){0.64018}}
  \put(7.9,0.35){\line(1,-3){0.31623}}
  \put(7.985,0.4){\line(2,-5){0.37139}}
  \put(8.07,0.42){\line(3,-5){0.5145}}
  \put(8.185,0.5){\line(4,-5){0.6247}}
  \put(8.27,0.55){\line(3,-4){0.6}}
  \put(8.47,0.7){\line(2,-3){0.5547}}
  \put(8.67,0.83){\line(3,-5){0.5145}} \thicklines
  \put(8.1,0.4){\circle{1.1}} \put(8.51,2.41){\vector(-1,-4){0.4}}
  \thinlines \put(4,3.3){there is no open neighborhood of $f(P(B_\delta^j))$, in}
  \put(4.5,2.8){which the normal projection is injective}}
\end{picture}
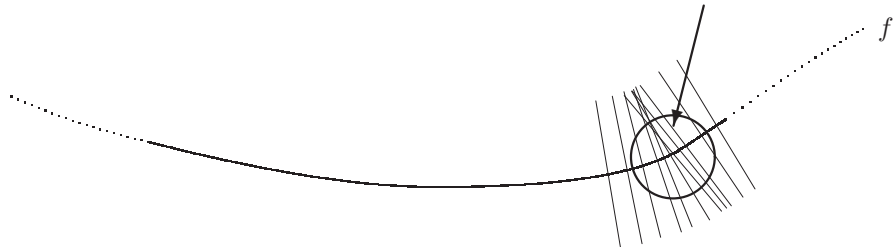
\\ \begin{fig} \label{normalconc} Normal projection in the case of concentrated curvature. \end{fig}
However, there are several ways for solving this problem. First one could smoothen the limit immersion
$f$ in order to obtain an immersion $g$, which is at least $C^2$ (or even $C^\infty$) and which is $C^1$-close
to $f$. Then we can project from $f$ in the normal direction $\nu_g$ of $g$ onto $f^i$. Similarly one could use
one of the approximation theorems for immersions in \cite{hirsch}. Slightly different is the approach using
an averaged normal projection. It is described in \cite{langer} for codimension $1$. A generalization
to arbitrary codimension using the Riemannian center of mass is presented in \cite{breuning1}. \\ \\
Here we like to assume that we have already found (by one of the preceding
methods) a smooth mapping $\nu:M\rightarrow G_{n,k}$, which is close to the normal
of $f$. Let us explain what that means: As explained above
$f(P(B_\delta^j))$ is the graph of a function $u$
on $B_\delta\subset \R^m\times \{0\}\subset\R^m\times \R^k$ with $\|Du\|_{C^0(B_\delta)}\leq \alpha$.
For $q\in P(B_\delta^j)$ consider the subspace $\nu_f(q)$, where $\nu_f:M\rightarrow G_{n,k}$ is the
normal of $f$. Then $\nu_f(q)$ is a graph over $\{0\}\times \R^k$;
more precisely there is a linear map $\tilde{N}_q:\R^k\rightarrow \R^m$ such that \beqnn
\nu_f(q)=\{(\tilde{N}_q(z),z):z\in\R^k\}\subset \R^m\times \R^k=\R^n. \eeqnn
Moreover, as $\|Du\|_{C^0(B_\delta)}\leq \alpha$, for the operator norm $\|\cdot\|_{\text{op}}$ we have \beqnn
\|\tilde{N}_q\|_{\text{op}}\:\leq\: \alpha. \eeqnn
The property of $\nu$ being close to $\nu_f$ (which can be reached by any of the
described methods) shall mean, that for
all $q\in P(B_\delta^j)$ also the subspace $\nu(q)$ is the graph of a linear map
$N_q:\R^k\!\rightarrow\! \R^m$ over $\{0\}\times \R^k$ and that
\beqn \|N_q\|_{\text{op}}\:\leq\: 2\alpha. \eeqn
Identifying $P(B_\delta^j)$ with the ball $B_\delta$ as described above,
we may similarly assign to each $x\in B_\delta$
a linear map $N_x:\R^k\rightarrow \R^m$. \\ \\
We like to show, that for $q \in P(B_\delta^j)$ the affine
subspace $f(q)+\nu(q)$ has exactly one point of intersection with
the set $f^i(U_{4\delta,j}^i)$. \\ \\
For that, in addition to (\ref{positionimmersion1}),
assume \beqn \label{graph2} \{(x,\tilde{u}^i(x)): x \in
B_{4\delta-\varrho}\}\subset f^i(U_{4\delta,j}^i), \eeqn where
$\varrho$ is small, say $\varrho=\frac{\delta}{2}$\: (suppose (\ref{positionimmersion1})
is satisfied with the same $\varrho$). \\ \\
\emph{The mapping $F$:} \\[1.5mm]
For $x \in B_\delta$ we denote by $F(x)$ the unique intersection point of
the affine subspace \beqnn h(x)\,:=\,(x,u(x))+\nu(x) \eeqnn with $\mathbb{R}^m\times
\{0\}$. In that way, we obtain the mapping \\ \parbox{12cm}{\beqnn F:B_\delta
&\rightarrow& \mathbb{R}^m,
\\ x &\mapsto& x-N_x(u(x)). \eeqnn} \hfill \parbox{8mm}{\beqn \eeqn}
\\ \\
\emph{The mappings $G_x^i$:} \\[1.5mm]
For $x \in B_\delta$ and $y \in B_{r-\varrho}$ we denote by
$G_x^i(y)$ the unique intersection point of the affine subspace
$(y,\tilde{u}^i(y))+\nu(x)$ with $\mathbb{R}^m \times \{0\}$. In
that way we obtain for each fixed $x \in B_\delta$\, a mapping \\
\parbox{12cm}{\beqnn G_x^i:B_{r-\varrho}&\rightarrow& \mathbb{R}^m, \\ y
&\mapsto&
y-N_x(\tilde{u}^i(y)). \eeqnn} \hfill \parbox{8mm}{\beqn \eeqn} \\[-5mm]
$\text{ }$ \hspace{6mm} \setlength{\unitlength}{1.2cm}
\begin{picture}(10,5)
\put(-0.6,0){\line(1,0){11.2}} \put(1.2,-0.065){(} \put(8.8,-0.065){)}
\put(-0.63,-0.065){(} \put(10.53,-0.065){)} \put(5,-0.4){0}
\put(8.7,-0.6){$B_\delta$} \put(10.4,-0.6){$B_{r-\varrho}$}
\put(10.8,2){$f$} \put(0.3,-0.065){(} \put(9.7,-0.065){)}
\put(9.4,-0.6){$B_{4\delta-\varrho}$} \put(-0.1,2.18){(}
\put(9.9,2.95){)}

\qbezier[20](-0.6,1.2)(0.3,0.8)(1.2,0.6)
     \qbezier(1.2,0.6)(3.5,0)(5,0) \qbezier(5,0)(7.4,0)(8.2,0.5)
\qbezier(8.2,0.5)(8.8,0.9)(8.8,0.9)
\qbezier[20](8.8,0.9)(9.9,1.7)(10.6,2.1)

\qbezier[300](-0.6,2.4)(6,0.7)(10.6,3.4) \put(10.8,3.4){$f^i$}

\put(1.7,0.45){\line(-2,-5){0.7}} \put(1.7,0.45){\line(2,5){1.3}}
\put(9.3,2.7){\line(-2,-5){1.6}} \put(9.3,2.7){\line(2,5){0.4}}

\put(9.3,-0.06){\tiny $|$} \put(1.7,-0.06){\tiny $|$}

\put(8.2,-0.06){\tiny $|$} \put(1.5,-0.06){\tiny $|$}
\put(9.335,-1.2){\vector(0,1){1}} \put(8.985,-1.5){$y \in
B_{r-\varrho}$} \put(8.235,-1.2){\vector(0,1){1}}
\put(7.85,-1.5){$G_x^i(y)$} \put(1.535,-1.2){\vector(0,1){1}}
\put(1.3,-1.5){$F(x)$} \put(2.8,-1.2){\vector(-1,1){1}}
\put(2.5,-1.5){$x \in B_\delta$}

\put(0.7,1){\vector(2,-1){1}} \put(0.2,1.15){$(x,u(x))$}
\put(8.25,3.25){\vector(2,-1){1}}
\put(7.65,3.4){$(y,\tilde{u}^i(y))$}

\put(2.5,4.5){\vector(1,-3){0.33}}
\put(0.7,4.7){$h(x)=(x,u(x))+\nu(x)$}
\put(9.2,4.5){\vector(1,-3){0.33}}
\put(7.9,4.7){$(y,\tilde{u}^i(y))+\nu(x)$}

\end{picture}
\vspace{1.8cm}
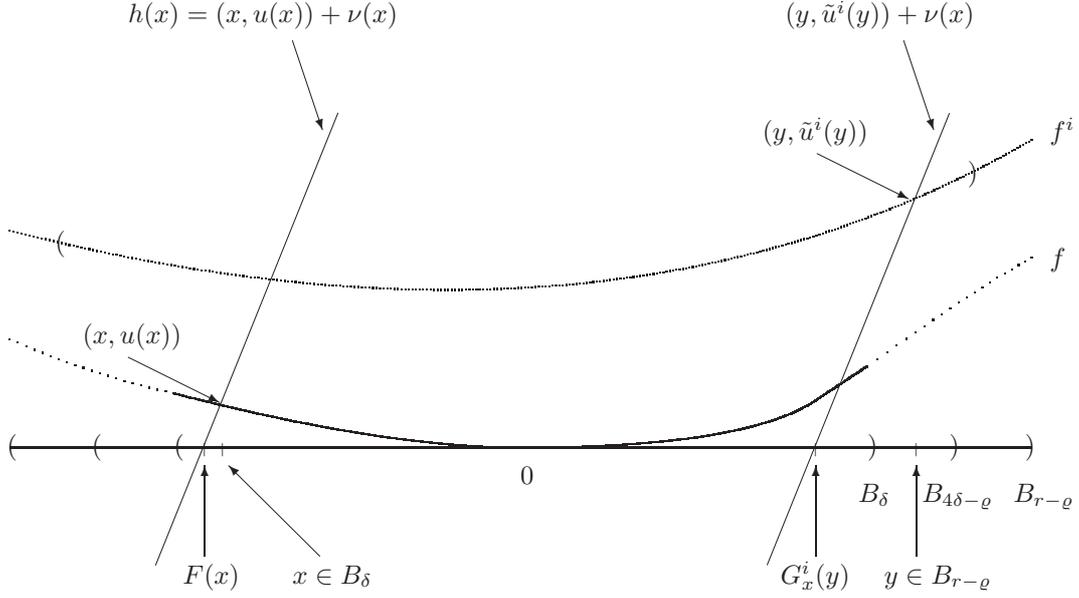
\begin{fig} The mappings $F$ and $G_x^i$. The part between the parentheses on the
immersion $f^i$ represents the set $f^i(U_{4\delta,j}^i)$.
\end{fig} \vspace{0cm}
\emph{The mappings $H_x^i$:} \\[1.5mm]
For $y \in \overline{B}_{4\delta-2\varrho}$ and $\varepsilon$ sufficiently small we have
$|\tilde{u}^i(y)|\leq 2\alpha(4\delta-2\varrho)+\varepsilon\leq 8\alpha\delta$, hence $|y-G_x^i(y)|
=|N_x(\tilde{u}^i(y))|\leq 16\alpha^2\delta$. For $x\in B_\delta$
we have $|u(x)|\leq \alpha\delta$, hence $|F(x)|=|x-N_x(u(x))|\leq \delta+2\alpha^2\delta$.
Using $\alpha^2\leq \frac{1}{10}$, $\varrho=\frac{\delta}{2}$, this yields \beqnn
|y-G_x^i(y)+F(x)|\;\leq\; (1+18\alpha^2)\delta \;\leq\; 3\delta\;=\;4\delta-2\varrho. \eeqnn
With that we define for each fixed $x\in B_\delta$\, a mapping \\
\parbox{12cm}{\beqnn H_x^i:\;\;
\overline{B}_{4\delta-2\varrho}&\rightarrow&
\overline{B}_{4\delta-2\varrho},
\\ y&\mapsto& y-G_x^i(y)+F(x). \eeqnn} \hfill \parbox{8mm}{\beqn \eeqn}
\begin{lemma}
\label{intersection} For $q \in P(B_\delta^j)$ the affine subspace
$f(q)+\nu(q)$ has exactly one point of intersection with the set
$f^i(U_{r,j}^i)$. This point lies in $f^i(U_{4\delta,j}^i)$.
\end{lemma}
\textbf{Proof:}\\
We pass to the local representation and consider $h(x)=(x,u(x))+\nu(x)$ for
$x\in B_\delta$.
Using the definition of $H_x^i$ and $\alpha^2\leq\frac{1}{10}$, we estimate \beqnn
|H_x^i(\xi)-H_x^i(\zeta)|&=&|N_x(\tilde{u}^i(\xi)-\tilde{u}^i(\zeta))| \\
&\leq& 2\alpha\|N_x\||\xi-\zeta| \\
&\leq& 4\alpha^2\,|\xi-\zeta| \\
&\leq& \frac{1}{2}\,|\xi-\zeta|.\eeqnn
Hence $H_x^i$ is a contraction. By the Banach fixed
point theorem there is exactly one \linebreak $y \in
\overline{B}_{4\delta-2\varrho}$ with $H_x^i(y)=y$, that is with
$G_x^i(y)=F(x)$. \\ \\
By the definitions of $F$ and $G_x^i$, the affine subspaces
$h(x)\!=(x,u(x))+\nu(x)$ and $(y,\tilde{u}^i(y))+\nu(x)$ intersect
each other in $F(x)=G_x^i(y)$ and are parallel, hence
$h(x)=(y,\tilde{u}^i(y))+\nu(x)$ and $(y,\tilde{u}^i(y))\in h(x)$.
By (\ref{graph2}), the affine subspace $h(x)$ intersects the set
$f^i(U_{4\delta,j}^i)$ in $(y,\tilde{u}^i(y))$. \\ \\[1mm]
Similarly, we show that there is only one point of intersection with $f^i(U_{r,j}^i)$:
For that we assume that we have chosen $r$ slightly
smaller in the beginning, such that also the set $f^i(U_{r+2\varrho,j}^i)$ is the
graph of a function $\tilde{u}^i$
on a subset of $\R^m$ with $\|D\tilde{u}^i\|_{C^0}\leq 2\alpha$, and such that
\beqnn
f^i(U_{r,j}^i)\:\subset\:
\{(x,\tilde{u}^i(x)): x \in \overline{B}_{r+\varrho}\,\}
\:\subset\:f^i(U_{r+2\varrho,j}^i).
\eeqnn
Now for each fixed $x\in B_\delta$ define a function
\\
\parbox{12cm}{\beqnn \widetilde{H}_x^i:\;\;
\overline{B}_{r+\varrho}&\rightarrow&
\overline{B}_{r+\varrho},
\\ y&\mapsto& y-G_x^i(y)+F(x), \eeqnn} \hfill
\\ \\
where we also extend $G_x^i$ to the ball $\overline{B}_{r+\varrho}$.
Using $r=16\delta$, $\varrho=\frac{\delta}{2}$, $\alpha^2\leq\frac{1}{10}$
and assuming $\varepsilon$ to be small, one shows
$|y-G_x^i(y)+F(x)|\leq \frac{r}{2}$. Hence $\widetilde{H}_x^i$ is well-defined. Then also
$\widetilde{H}_x^i$ is a contraction and there \nolinebreak is exactly one
$y\in\overline{B}_{r+\varrho}$ with $G_x^i(y)=F(x)$. By the definitions of
$G_x^i$ and $F$,
this shows the statement. \hfill $\square$ \\ \\
Before we come to the definition of the mappings $\phi^i:M\rightarrow M^i$, we need
the following lemma, which will assure that the $\phi^i$ are well-defined:
\begin{lemma} \label{phiwell}
Let $x\in P(B_\delta^j)\cap P(B_\delta^k)$. Moreover let $S_1$ be the point of intersection
of $h(x)$ with $f^i(U_{r,j}^i)$, $S_2$ the point of intersection of $h(x)$ with
$f^i(U_{r,k}^i)$, and $\sigma_1 \in U_{r,j}^i$ with $f^i(\sigma_1)=S_1$,
$\sigma_2\in U_{r,k}^i$ with $f^i(\sigma_2)=S_2$. Then $\sigma_1=\sigma_2$.
\end{lemma}
\textbf{Proof:} \\
By Lemma \ref{intersection} we have $S_2\in f^i(U_{4\delta,k}^i)$, that is
$\sigma_2\in U_{4\delta,k}^i$. The assumption $x\in P(B_\delta^j)\cap P(B_\delta^k)$
implies $k\in Z(j)$, hence by Lemma \ref{inclusiondipl} b)
$U_{4\delta,k}^i\subset U_{r,j}^i$. Using again Lemma \ref{intersection},
the statement follows. \hfill $\square$ \\ \\
With the preceding lemmas we are able to give a definition of the mappings
$\phi^i:M\rightarrow M^i$. For that let $x\in M$.  Then $x\in P(B_\delta^j)$
for some $j$. The set $h(x)$ intersects $f^i(U_{r,j}^i)$ in exactly one
point $S_x$. Furthermore there is exactly one point $\sigma_x\in U_{r,j}^i$
with $f^i(\sigma_x)=S_x$.
We set $\phi^i(x):=\sigma_x$. The mappings $\phi^i$ are well-defined by
Lemma \ref{phiwell}.
Now we like to show that the mappings $\phi^i$ (after passing to a subsequence,
if necessary) are diffeomorphisms. \\ \\
Let $\gamma_{n,k}=\{(E,x):E\in G_{n,k}, x\in E\}$ and
let $p:\gamma_{n,k}\rightarrow G_{n,k}$, $(E,x)\mapsto E$,
be the universal bundle over $G_{n,k}$. The local trivializations for this bundle are defined as follows:
Let $E\in G_{n,k}$ and let $\pi_E:\R^n\rightarrow E$ be the orthogonal projection; we set
$U_E=\{G\in G_{n,k}:\pi_E(G) \text{ is of dimension } k\}$; \linebreak a local trivialization
is then given by $\Psi:p^{-1}(U_E)\rightarrow U_E\times E\cong U_E\times \R^k$,
$\Psi((G,x))=(G,\pi_E(x))$. \\ \\
Let $\nu:M\rightarrow G_{n,k}$ be as above.
We now consider the pullback bundle $\nu^\ast \gamma_{n,k}$, which is a vector bundle over $M$ with bundle
projection $\pi$ and $n$-dimensional total space \beqnn
E=\{(x,y)\in M\times \R^n:y\in \nu(x)\}. \eeqnn
We set $E_j=\{(x,y)\in E:x \in P(B_\delta^j)\}$. Hence $\nu^\ast \gamma_{n,k}|E_j$ is a bundle
over $P(B_\delta^j)$. As $P(B_\delta^j)$ is diffeomorphic to $B_\delta\subset \R^m$,
and as $B_\delta$ is diffeomorphic to $f(P(B_\delta^j))=\{(x,u_j(x)):x\in B_\delta\}$,
we may consider $\nu^\ast \gamma_{n,k}|E_j$ also as a bundle over one of the last-named
sets. In particular, $\nu^\ast \gamma_{n,k}|E_j$ is a trivial bundle. \\ \\
We sometimes identify the zero section of $\nu^\ast \gamma_{n,k}|E_j$ with $P(B_\delta^j)$.
Finally we define a mapping
\\[-3mm] \parbox{12cm}{\beqnn F:\;\;\;\;\; \,E\;\;&\rightarrow& \mathbb{R}^n,\\
(x,y)&\mapsto&f(x)+y, \eeqnn} \hfill
\parbox{8mm}{\beqn \eeqn} \\
where $y \in \nu(x)$.

\begin{lemma}[Local tubular neighborhood around the limit immersion]
\label{tubular}
There exists an open neighborhood $V\subset E$
of the zero section of $\nu^{\ast}\gamma _{n,k}$, such that
for every $j$ with $1\leq j \leq s$ the following holds:\begin{itemize}
\item $F|E_j\cap V$ is a diffeomorphism onto an open neighborhood of $f(P(B_\delta^j))$,
\item $F|P(B_\delta^j)=f|P(B_\delta^j)$,
\item for every fibre $E_q=\pi^{-1}(q)$ we have $F(E_q)=h(q)$.
\end{itemize}
\end{lemma}
\textbf{Proof:} \\
We note that for every $q\in M$ the affine subspace $f(q)+\nu(q)$
intersects $f(q)$ transversally. Moreover $\nu$ is a smooth mapping.
Now the statement is a simple fact from differential topology about the
existence of tubular neighborhoods (see \cite{broecker} and \cite{hirsch}). In this
way we find tubular neighborhoods on $P(B_\delta^j)$ for every $j$.
Appropriately composing these neighborhoods, we obtain the desired
neighborhood $V\subset E$ of the zero section of $\nu^{\ast}\gamma _{n,k}$. \hfill $\square$ \\[-3mm]
\begin{lemma} \label{surj} After passing to a subsequence, each mapping $\phi^i:M\rightarrow M^i$
is surjective. \end{lemma}
\textbf{Proof:} \\
By Lemma \ref{tubular}, for each $j$ the set $F(E_j\cap V)$ is an open neighborhood
of $f(P(B_\delta^j))=\{(x,u_j(x)):x \in B_\delta\}$. We define sets
$M_j=\{(x,u_j(x)):x \in \overline{B_{\frac{2}{3}\delta}}\}\subset
f(P(B_\delta^j))$. As $M_j$ is compact, there is an $\varepsilon_j>0$
with \vspace{-3mm} \beqnn M_j^{\varepsilon_{j}} :=\{(x,y): x \in
\overline{B_{\frac{2}{3}\delta}}, \; \; y\in \mathbb{R}^k \text{ mit
} |y-u_j(x)|<\varepsilon_{j}\}\subset F(E_j\cap V). \eeqnn We set
$\hat{\varepsilon}=\min\{\varepsilon_1,\ldots,\varepsilon_s\}$. By
definition of $\mathfrak{d}(\cdot,\cdot)$
and by graph convergence, it follows that
$f^i(U_{\delta/2,j}^{i})$ is a subset of
$M_j^{\hat{\varepsilon}}$ for $i$ sufficiently large
(see Figure \ref{mappingsurjective}). Further, it follows
$U_{\delta/2,j}^{i}\subset \phi^i(P(B_\delta^j))$ for
$j=1,\ldots,s$. Hence, for every $q\in U_{\delta/2,j}^i$ there is a
$p\in P(B_\delta^j)$ with $f^i(q)\in f(p)+\nu(p)$. By the definition
of $\phi^i$, this yields $\phi^i(p)=q$.
As $\{q_1^i,\ldots,q_s^i\}$ is a $\frac{\delta}{2}$-net for $f^i$, for every $q\in M^i$ there is a $j\in\{1,\ldots,s\}$
with $q\in U_{\delta/2,j}^i$. Hence, by the considerations of above, $\phi^i$ is surjective. \hfill $\square$ \\ \vspace{0.5cm} \\
\setlength{\unitlength}{1.2cm}
\begin{picture}(10,3.2) \put(1.5,0){
\put(10.8,2){$f$}

\put(-0.6,0){\line(1,0){11.2}} \put(1.2,-0.1){(} \put(8.8,-0.1){)}
\put(-0.6,-0.1){(} \put(10.5,-0.1){)} \put(5,-0.4){0}
\put(8.7,-0.6){$B_\delta$} \put(10.4,-0.6){$B_r$}
\put(2.45,-0.1){\textbf{[}} \put(7.47,-0.1){\textbf{]}}
\put(2.2,0.7){$\bigg \{$} \put(2,0.7){$\hat{\varepsilon}$}
\put(8.15,-1.5){$\overline{B_{\frac{2}{3}\delta}}$}
\put(8.3,-1.1){\vector(-2,3){0.7}}
\put(2.9,2.5){$f^i(U_{\frac{\delta}{2},j}^i)$}

\put(3.3,2.27){\vector(1,-3){0.55}}

\qbezier[20](-0.6,1.2)(0.3,0.8)(1.2,0.6)
     \qbezier(1.2,0.6)(3.5,0)(5,0) \qbezier(5,0)(7.4,0)(8.2,0.5)
\qbezier(8.2,0.5)(8.8,0.9)(8.8,0.9)
\qbezier[20](8.8,0.9)(9.9,1.7)(10.6,2.1)

\put(8.85,0.85){\line(2,-5){1.11417}}
\put(8.85,0.85){\line(-2,5){0.95}}

\put(1.2,0.55){\line(-1,-4){0.72762}}
\put(1.2,0.55){\line(1,4){0.72762}}

\linethickness{0.5mm} \qbezier(1.3,1.1)(6.2,2.6)(8.65,1.3)
\qbezier(0.75,-1.2)(5.1,-4)(9,0.5) \thinlines
\qbezier(2.467,1.3)(5.1,0.68)(7.534,1.25)

\qbezier(2.467,-0.7)(5.1,-1.32)(7.534,-0.75)

\qbezier(2.8,0.8)(4.5,0.2)(6.6,0.5)

\thicklines \put(2.467,-3.4){\line(0,1){7}}
\put(7.534,-3.4){\line(0,1){7}}}
\end{picture}
\vspace{3.8cm}
\begin{fig} \label{mappingsurjective} Surjectivity of the mappings $\phi^i$. \end{fig} \pagebreak
For showing injectivity, we need the following lemma:
\begin{lemma}\label{help}
For $i$ sufficiently large, we have the inclusions \\
a) $\phi^i(P(B_{\frac{\delta}{3}}^j))\subset
U_{\frac{\delta}{2},j}^i$, \\
b) $\phi^i(P(B_{\frac{\delta}{6}}^j))\subset
U_{\frac{\delta}{5},j}^i$.
\end{lemma}
\textbf{Proof:} \\
Follow the arguments of Lemma \ref{intersection}. A calculation with the numbers of above
proves a) and b). \hfill $\square$ \\ \\
We first show local injectivity:
\begin{lemma}[Local injectivity] \label{locinj}
After passing to a subsequence, for each $j$ the mappings $\phi^i:M\rightarrow M^i$ restricted to
$P(B_\delta^j)$ are injective. \end{lemma}
\textbf{Proof:} \\
Let $x,y \in P(B_\delta^j)$ with $x \neq y$. By the definition of
$\phi^i$ and by Lemma \ref{tubular} we have
$f^i\circ \phi^i(x) \in h(x)=F(E_x)$, $f^i\circ \phi^i(y) \in
h(y)=F(E_y)$, and \beqn \label{leereMenge} F(E_x\cap V)\cap
F(E_y\cap V)=\emptyset. \eeqn By Lemma \ref{help} a) we have
$\phi^i(P(B_{\delta/3}^k))\subset U_{\delta/2,k}^i$ for each
$k$, which implies for $z \in P(B_{\delta/3}^k)$ with the arguments
in the proof of Lemma \ref{surj}, that
$f^i\circ \phi^i(z) \in F(E_z\cap V)$. As by Lemma
\ref{eigenschaft1} it holds $M \subset
\bigcup_{j=1}^{s}P(B_{\delta/3}^k)$, we actually have
$f^i\circ \phi^i(z)\in F(E_z\cap V)$ for all $z \in M$. With
(\ref{leereMenge}) it follows $f^i\circ \phi^i(x) \neq f^i\circ
\phi^i(y)$, and hence $\phi^i(x) \neq \phi^i(y)$. \hfill $\square$ \\ \\
Now we like to show global injectivity:
\begin{lemma}[Injectivity of \boldmath $\phi^i$
\unboldmath]\label{inj} After passing to a subsequence, the
mappings $\phi^i:M\rightarrow M^i$ are injective.
\end{lemma}
\textbf{Proof:} \\
Let $x,y \in M$ with $x \neq y$. By Lemma \ref{eigenschaft1} there are
$j,k$ with $x \in P(B_{\delta/6}^j)\subset P(B_{\delta/4}^j)$, $y
\in
P(B_{\delta/6}^k)\subset P(B_{\delta/4}^k)$. \\ \\
\textbf{Case 1:} $P(B_{\frac{\delta}{4}}^j)\cap
P(B_{\frac{\delta}{4}}^k) =
\emptyset$ \\ \\
By Lemma \ref{help} b) we have $\phi^i(x) \in U_{\delta/5,j}^i$,
$\phi^i(y) \in U_{\delta/5,k}^i$ and Lemma \ref{eigenschaft3}
implies $k \notin \tilde{Z}(j)$, that is $U_{\delta/5,j}^i \cap
U_{\delta/5,k}^i =
\emptyset$. It follows $\phi^i(x) \neq \phi^i(y)$. \\ \\
\textbf{Case 2:} $P(B_{\frac{\delta}{4}}^j)\cap
P(B_{\frac{\delta}{4}}^k)
\neq \emptyset$ \\ \\
By Lemma \ref{eigenschaft2} we have $P(B_{\delta/4}^k)\subset
P(B_{\delta}^j)$. By Lemma \ref{locinj} $\phi^i$ is injective
on $P(B_{\delta}^j)$, hence again $\phi^i(x) \neq \phi^i(y)$. \hfill $\square$ \\ \\
For showing, that each mapping $\phi^i$ is a diffeomorphism, we first show that the
composition $f^i\circ \phi^i$ is an immersion. For that, we use that $F(E_j\cap V)$
is a tubular neighborhood both of $f(P(B_\delta^j))$ and of $f^i\circ \phi^i(P(B_\delta^j))$.
\begin{lemma} \label{compimmersion} The mapping $f^i\circ \phi^i:M\rightarrow \R^n$ is an immersion. \end{lemma}
\textbf{Proof:} \\
We show the statement by considering the local representation of
$f^i\circ \phi^i$ on the set $P(B_\delta^j)$. We
regard $E_j$ as a bundle over $f(P(B_\delta^j))$. As
$\nu_{g}^{*}\gamma_{n,k}|E_{j}$ is a trivial bundle,
there exists a trivialization $\Psi_j:E_j\rightarrow
B_{\delta}\times \mathbb{R}^k \subset \mathbb{R}^n$ with
$\Psi_j(E_q)=\{q\}\times \mathbb{R}^k$. As $f(P(B_\delta^j))$
is diffeomorphic to $B_\delta$, we may assume that the zero section
is mapped by $\Psi_j$ onto $B_\delta$. We
define restrictions $\widetilde{\Psi}_j=\Psi_j|E_j\cap V:
\;E_j\cap V\rightarrow \Psi_j(E_j\cap V)$, and $F_j=
\;F|E_j\cap V: \;\; E_j\cap V\rightarrow F(E_j\cap V)$ with $F$
as in Lemma \ref{tubular}. Also by Lemma
\ref{tubular}, $F_j$ and hence also
$\widetilde{\Psi}_j\circ F_{j}^{-1}$ is a
diffeomorphism (see Figure \ref{geradebiegen}).
We note, that $f^i(U_{r,j}^i)$ and hence also
$W_j^i:=f^i(U_{r,j}^i)\cap F(E_j\cap V)$ is a smooth submanifold
of $\mathbb{R}^n$.
With Lemma \ref{tubular} and
by construction of the projection, for $q \in P(B_\delta^j)$ we have $f^i\circ \phi^i(q) \in
f^i(U_{r,j}^{i})\cap h(q)=f^i(U_{r,j}^{i})\cap F(E_q)$, and $f^i\circ
\phi^i(B_\delta)=W_j^i$.
We obtain $F_j^{-1}\circ f^i\circ \phi^i(q) \in E_q$ and hence
$\widetilde{\Psi}_j \circ F_j^{-1}\circ f^i \circ
\phi^i(q)=(q,h_j^i(q)) \in B_\delta \times \mathbb{R}^k$ with a
mapping
$h_j^i:B_\delta\rightarrow \mathbb{R}^k$.
As $\widetilde{\Psi}_j \circ F_j^{-1}$ is a diffeomorphism and $W_j^i$
a submanifold, also $\widetilde{\Psi}_j \circ
F_j^{-1}(W_j^i) =\widetilde{\Psi}_j \circ F_j^{-1}\circ
f^i\circ \phi^i(B_\delta)$ is a smooth submanifold and hence
$h_j^i$ a differentiable mapping. It follows, that
$\widetilde{\Psi}_j \circ F_j^{-1}\circ f^i\circ \phi^i$ is an immersion.
As $\widetilde{\Psi}_j\circ F_j^{-1}$ is a diffeomorphism, the statement follows.
\hfill $\square$ \\  $\text{ }$ \hspace{0.5cm}
\setlength{\unitlength}{0.65cm}
\begin{picture}(15,5.5) \put(1.5,0){

\qbezier[20](-2,1)(-1.7,0.8)(-0.8,0.6)
     \qbezier(-0.8,0.6)(1.5,0)(3,0) \qbezier(3,0)(5.4,0)(6.2,0.5)
\qbezier(6.2,0.5)(6.8,0.9)(6.8,0.9)
\qbezier[20](6.8,0.9)(7.3,1.25)(7.5,1.45)

\put(6.85,0.85){\line(2,-5){1.11417}}
\put(6.85,0.85){\line(-2,5){0.95}}

\put(6.35,0.55){\line(2,-3){1.6641}}\put(6.35,0.55){\line(-2,3){1.6641}}

\put(5.85,0.3){\line(2,-5){1.11417}}
\put(5.85,0.3){\line(-2,5){1.11417}}

\put(5.35,0.15){\line(1,-2){1.34163}}
\put(5.35,0.15){\line(-1,2){1.34163}}

\put(4.85,0.05){\line(2,-5){1.11417}}
\put(4.85,0.05){\line(-2,5){1.11417}}

\put(4.35,0){\line(1,-3){0.94869}}
\put(4.35,0){\line(-1,3){0.94869}}

\put(3.85,0){\line(1,-6){0.4932}} \put(3.85,0){\line(-1,6){0.4932}}

\put(3.35,-0.05){\line(1,-5){0.58836}}
\put(3.35,-0.05){\line(-1,5){0.58836}}

\put(2.85,-0.05){\line(1,-6){0.4932}}
\put(2.85,-0.05){\line(-1,6){0.4932}}

\put(2.3,0){\line(0,-1){3}} \put(2.3,0){\line(0,1){3}}

\put(1.795,0.05){\line(-1,-5){0.58836}}
\put(1.795,0.05){\line(1,5){0.58836}}

\put(1.3,0.1){\line(-1,-6){0.4932}}
\put(1.3,0.1){\line(1,6){0.4932}}

\put(0.83,0.2){\line(-1,-4){0.72762}}
\put(0.83,0.2){\line(1,4){0.72762}}

\put(0.3,0.28){\line(-1,-5){0.58836}}
\put(0.3,0.28){\line(1,5){0.58836}}

\put(-0.25,0.37){\line(-1,-6){0.4932}}
\put(-0.25,0.37){\line(1,6){0.4932}}

\put(-0.8,0.55){\line(-1,-4){0.72762}}
\put(-0.8,0.55){\line(1,4){0.72762}}

\qbezier(-2,1.2)(5.5,-0.4)(7.5,1.7)

\linethickness{0.5mm} \qbezier(-0.7,1.1)(4.2,2.6)(6.65,1.3)
\qbezier(-1.25,-1.2)(3.1,-4)(7,0.5)

\thinlines \put(12,0.3){\line(1,0){7.6}}

\put(12,-2.7){\line(0,1){6}} \put(12.507,-2.7){\line(0,1){6}}
\put(13.013,-2.7){\line(0,1){6}} \put(13.52,-2.7){\line(0,1){6}}
\put(14.027,-2.7){\line(0,1){6}} \put(14.533,-2.7){\line(0,1){6}}
\put(15.04,-2.7){\line(0,1){6}} \put(15.547,-2.7){\line(0,1){6}}
\put(16.053,-2.7){\line(0,1){6}} \put(16.56,-2.7){\line(0,1){6}}
\put(17.067,-2.7){\line(0,1){6}} \put(17.573,-2.7){\line(0,1){6}}
\put(18.08,-2.7){\line(0,1){6}} \put(18.587,-2.7){\line(0,1){6}}
\put(19.093,-2.7){\line(0,1){6}} \put(19.6,-2.7){\line(0,1){6}}

\linethickness{0.5mm} \qbezier(12,1.1)(15,2.8)(19.6,1.3)
\qbezier(12,-1.3)(14.9,-3.2)(19.6,-0.7)

\thinlines \qbezier(12,0.8)(15,1.3)(19.6,0.8)

\qbezier(7.6,2.5)(9.6,4)(11.2,2.6)
\put(11.2,2.6){\vector(1,-1){0.1}}

\put(8.5,3.55){$\widetilde{\Psi}_j\circ F_j^{-1}$}

\put(-2.1,-0.5){\vector(3,4){0.8}} \put(-2.45,-0.75){$f$}

\put(-2.1,2.2){\vector(3,-4){0.8}} \put(-2.45,2.2){$f^i$}

\put(14.8,3.8){$B_\delta \times \mathbb{R}^k$}}

\end{picture}

\vspace{1.9cm}
\begin{fig} \label{geradebiegen}
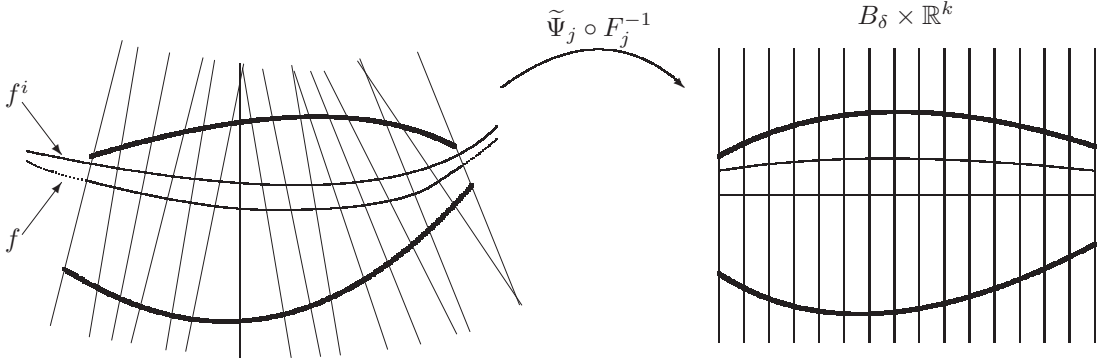
Straightening of the tubular neighborhood. Here the set $f(P(B_\delta^j))$ is mapped by
$\widetilde{\Psi}_j\circ F_j^{-1}$ onto $B_\delta\times \{0\}$.
Note that for $\nu \in C^{2}$ the mapping $f^i\circ \phi^i$ is in $W^{2,p}$.\end{fig}
\begin{theorem} The mappings $\phi^i:M\rightarrow M^i$ are diffeomorphisms. \end{theorem}
\textbf{Proof:} \\
The mappings $f^i$ and $f^i\circ \phi^i$ are immersions. It follows, that also $\phi^i$ is an immersion.
Moreover $\phi^i$ is surjective by Lemma \ref{surj} and injective by Lemma
\ref{inj}. Hence $\phi^i$ is a diffeomorphism. \hfill $\square$
\end{section}
\begin{section}{Convergence of the immersions} \label{convergenceimm}
In this section we would like to show convergence of the sequence
$f^i\circ \phi^i$ to $f$ in the $C^1$-topology. This means, we show
$C^1$-convergence of the local representations of $f^i\circ\phi^i$
to the local representations of the limit immersion $f$ with respect to
the atlas $\mathfrak{A}$. \\ \\
As a generalization, we like to show higher order convergence for immersions with
graph representations that are uniformly bounded in $W^{\K,p}$ with $\K>2$.
For that reason, let us assume that the mappings $f$ and $f^i$, and hence
also $u_j$, $u_j^i$ and $\tilde{u}_j^i$ are in $W^{\K,p}$ for
a $\K\geq 2$, and that $\tilde{u}_j^i$ is uniformly bounded in $W^{\K,p}$. We will
discuss in the end of this section under which assumptions we obtain these higher order bounds. \\ \\
We shall use the same notation as in the previous section. All considerations
are performed locally on $P(B_\delta^j)$. As in the previous section, we will fix
the index $j$ and suppress it in the notation. Hence again we shall write
$u$ instead of $u_j$ and $\tilde{u}^i$ instead of $\tilde{u}_j^i$. \\ \\
Instead let a lower index $\nu$ now denote the $\nu$-th coordinate of a vector
in $\R^n$ or in $\R^k$. Moreover, let $\pi^h$ be the projection from $\R^n$
onto the first $m$ coordinates, and $\pi^v$ the projection onto the last
$k$ coordinates. Finally, we shall simply write $f(x)$ instead of $(x,u(x))$. \\ \\
For proving convergence, we additionally assume that we have chosen
\beqn \label{alpha1k}
\alpha\leq \frac{1}{4\sqrt{k}} \eeqn
in the beginning, where $k$ denotes the codimension of the immersions
(here we denote by $\K$ the degree of differentiability,
and by $k$ the codimension). \\ \\
By the previous section, we project into the direction $\nu$. Moreover
$\nu^{*}\gamma_{n,k}|E_j$ is a trivial bundle over $B_\delta$. The fibre
of this bundle over each point in $B_\delta$ is a $k$-dimensional
subspace of $\R^n$. Now let
$\sigma=(e^1,\ldots,e^k)$ be a smooth frame of this bundle,
that is $e^1,\ldots,e^k:B_\delta \rightarrow \mathbb{R}^n$
and $(e^1(x),\ldots,e^k(x))$ is a basis of $\nu(x)\in
G_{n,k}$ for all $x \in B_\delta$. \\ \\
We define a mapping \beqn
G^i(x,t_1,\ldots,t_k)=\tilde{u}^i(\pi^h(f(x)+\sum
\limits_{\nu=1}^{k}t_{\nu}e^\nu))-\pi^v(f(x)+\sum
\limits_{\nu=1}^{k}t_\nu e^\nu), \eeqn where $x\in B_\delta$
and $t_1,\ldots,t_k \in \mathbb{R}$ is sufficiently small, such that
$\pi^h(f(x)+\sum \limits_{\nu=1}^{k}t_{\nu}e^\nu) \in B_{r-\varrho}$. \\ \\
By the construction of the reparametrization in the previous section,
for every $x \in B_\delta$ there is exactly one tuple
$(T_1^i(x),\ldots,T_k^i(x))$ with \beqn
G^i(x,T_1^i(x),\ldots,T_k^i(x))=0. \eeqn
In this manner we obtain mappings
$T_\nu^i:B_\delta \rightarrow \mathbb{R}$
(depending on the choice of frame). We like to choose a frame, such that
all calculations get as simple as possible. \\ \\
For that let $\hat{e}^1,\ldots,\hat{e}^k$ denote the standard orthonormal basis of
$\{0\}\times \mathbb{R}^k\subset
\mathbb{R}^m\times \mathbb{R}^k=\mathbb{R}^n$. For every $x
\in B_\delta$ the $k$-space $\nu(x)$ is a graph over
$\{0\}\times\mathbb{R}^k$. Now we define a frame
$\sigma=(e^1,\ldots,e^k)$ for
$\nu^{*}\gamma_{n,k}|E_j$ by projecting the basis
$\hat{e}^1,\ldots,\hat{e}^k$ orthogonally with respect to
$\{0\}\times\mathbb{R}^k$ onto each fibre $\nu(x)\in G_{n,k}$.
If $\nu$ is a $C^{\mathrm{k}}$-mapping, the basis vectors
$e^j:B_\delta\rightarrow \mathbb{R}^n$
are easily seen to be mappings of the same class. Moreover (as the bundle
$\nu^\ast \gamma_{n,k}|E_j$ can be continued to a trivial bundle on a larger set),
the mappings $e^j$ are bounded in $C^{\mathrm{k}}$. \\ \\
By construction, this frame has the property
$\pi^v(\sum_{\nu=1}^{k}t_\nu e^\nu)=(t_1,\ldots,t_k)^t$. For the corresponding
mappings $T_1^i,\ldots,T_k^i$, using
$\pi^h(f(x))=x$, it follows \beqn \label{vorkoordinate}
0=\tilde{u}^i(Id_{B_{\delta}}+\pi^h(\sum \limits_{\nu=1}^{k}T_\nu^i
e^\nu))-u-(T_1^i,\ldots,T_k^i)^t, \eeqn
that is for each coordinate
\beqn \label{koordinate}
0=\tilde{u}_\iota^i(Id_{B_{\delta}}+\pi^h(\sum
\limits_{\nu=1}^{k}T_\nu^i e^\nu))-u_\iota-T_\iota^i \eeqn for
all $\iota$ with $1\leq
\iota \leq k$. \\ \\
Now $\tilde{u}^i(Id_{B_{\delta}}+\pi^h(\sum _{\nu=1}^{k}T_\nu^i
e^\nu))$ is just the local representation of the mapping
$f^i\circ \phi^i$, which is in $W^{2,p}$. It follows directly,
that also the mappings $T_\iota^i$ are in $W^{2,p}$.
As we also like to show higher order convergence, let us assume that each
mapping $T_\iota^i$ is in $W^{\mathrm{k},p}$ with $\mathrm{k} \geq 2$. \pagebreak \\
\textbf{Convergence} \\
We like to show, that the mappings $T_\iota^i$ are uniformly bounded
in $W^{\mathrm{k},p}$. \\ \\
Beforehand we define \beqn X^i:B_\delta\rightarrow
B_{r-\varrho}, \;\;\;\;\;\;\;X^i= Id_{B_{\delta}}+\pi^h(\sum
_{\nu=1}^{k}T_\nu^i e^\nu). \eeqn Here all balls are subsets of
$\mathbb{R}^m$. Inserting $X^i$ into
(\ref{koordinate}) gives
\beqn \label{xieingesetzt} T_\iota^i=\tilde{u}_\iota^i\circ
X^i-u_\iota. \eeqn
\\
The following expressions $S_l^i$ and $U_l^i$ need not be calculated explicitly;
we only need to estimated the order of derivatives involved.
For that we shall use the multi-index notation. Expressions of the form
$\partial_w$ will denote the usual partial derivative.
Lower indices (as, for example, in the expressions
$X_\alpha^i$, $T_\nu^i$, $e_\alpha^\nu$, $u_\iota$,
$\tilde{u}_\iota^i$) will
denote the corresponding coordinate of a vector.
\begin{lemma} \label{ableitung1} Let $\gamma \in \mathbb{N}_0^m$, $|\gamma|=l$, be a multi-index with
$1\leq l \leq \mathrm{k}$. Then \beqn \label{ablgl1}
\partial^{\gamma}X_{\alpha}^{i}=\sum
\limits_{\nu=1}^{k}\partial^{\gamma}T_{\nu}^{i}\cdot
e_{\alpha}^{\nu}+S_l^i, \eeqn where $S_l^i$ is a finite sum of terms of the form
$C+\partial^{\lambda}T_\nu^i\cdot
\partial^{\mu}e_{\alpha}^{\nu}$ with multi-indices $\lambda,\mu$
with $0\leq|\lambda|\leq l-1$, \hspace{1.5mm} $1\leq |\mu| \leq l$
and $C$ a constant.
\end{lemma}
\textbf{Proof:} \\
The statement is easily shown by induction over $l$ (the order of the multi-index
$\gamma$). \hfill $\square$ \\
\begin{lemma} \label{ableitung2} Let $\gamma \in \mathbb{N}_0^m$, $|\gamma|=l$, be a multi-index with
$1\leq l \leq \mathrm{k}$. Then \beqn \label{ablgl2}
\partial^\gamma T_\iota^i = \sum
\limits_{\alpha=1}^{m}\partial_{\alpha}\tilde{u}_\iota^i(X^i)\cdot\partial^{\gamma}X_\alpha^i
-\partial^{\gamma}u_\iota+U_l^i , \eeqn where $U_l^i$ is a finite sum of terms of the form
$\partial^\lambda \tilde{u}_\iota^i(X^i)\cdot
\partial^{\mu_{1}}X_{\beta_{1}}^{i}\cdot \ldots \cdot
\partial^{\mu_{\eta}}X_{\beta_{\eta}}^{i}$ with $1\leq|\lambda|
\leq l$, \hspace{1.5mm} $1 \leq \eta \leq l$, \hspace{1.5mm} $1\leq
|\mu_1|,\ldots,|\mu_\eta| \leq l-1$ and $1 \leq \beta_1,\ldots,
\beta_\eta \leq m$.
\end{lemma}
\textbf{Proof:} \\
Again the statement is shown by induction over $l$. For $l=1$ and
 $1 \leq w \leq m$ one calculates the derivative
$\partial_w$ of equation (\ref{xieingesetzt}). The induction step is shown by straightforward
calculations. \hfill $\square$ \\ \\
Before showing convergence in $C^{\K-1}$, we show pointwise convergence:
\begin{lemma} \label{punktweise} It holds pointwisely $T^i\rightarrow 0$ as $i\rightarrow
\infty$. \end{lemma}
\textbf{Proof:} \\
Let $x\in B_\delta$ and $\varepsilon > 0$. By the graph convergence
there is an $N \in \mathbb{N}$, such that
 \beqn \label{konvpkt}
 \|\tilde{u}^i-u\|_{C^{0}(B_{r-\varrho})}<\frac{\varepsilon}{2}
 \;\;\;\;\;\;
\text{ for all } \;\; i>N, \eeqn
where $u$ is the corresponding function of the limit graph system $\Gamma$.
Let $y^i$ be the local representation of the point
$f^i\circ \phi^i(x)$, that is
$y^i=f(x)+\sum_{\nu=1}^{k}T_\nu^i(x)e^\nu(x)$, $y^i=(y_i^h,y_i^v)\in
\mathbb{R}^m\times \mathbb{R}^k$. By construction of the mappings
$\phi^i$, we have $y^i=(y_i^h,\tilde{u}^i(y_i^h))$. We set \beqnn
\varepsilon^i=|f^i\circ \phi^i(x)-f(x)|. \eeqnn The slope of
$\nu_g(x)$, and $\|Du\|_{C^{0}(B_{r})}\leq \alpha$ imply
$|\tilde{u}^i(y_i^h)-u(y_i^h)|>\frac{\varepsilon^i}{2}$. With
(\ref{konvpkt}) it follows $\varepsilon^i<\varepsilon$ for all $i>N$.\\ \\
Hence
$\varepsilon^i=|\sum_{\nu=1}^{k}T_\nu^i(x)e^\nu(x)|\rightarrow 0$
as $i\rightarrow \infty$. As the vectors
$e^1(x),\ldots,e^k(x)$ are linearly independent, we finally conclude
$T^i(x)\rightarrow 0$ as $i \rightarrow
\infty$. \hfill $\square$ \\ \\
Now we are able to show convergence in $C^{\K-1}$:
\begin{theorem} \label{konvergenz}
Under the assumptions at the beginning of this section, a subsequence
of $f^i\circ \phi^i$ converges in the $C^{\mathrm{k}-1}$-topology to $f$.
In particular, it follows $C^1$-convergence in the situation of Theorem
\ref{compactness0}.
\end{theorem}
\textbf{Proof:} \\
Let $\gamma$ be as in Lemmas \ref{ableitung1} and \ref{ableitung2}. We
insert (\ref{ablgl1}) into (\ref{ablgl2}) and obtain
\\ \parbox{12cm}{\beqnn \hspace{2.33cm} \partial^\gamma
T_\iota^i&=&\sum_{\alpha=1}^m\partial_\alpha\tilde{u}_\iota^i(X^i)\cdot\left(\sum_{\nu=1}^k
\partial^\gamma T_\nu^i \cdot e_\alpha^\nu+S_l^i\right)-\partial^\gamma
u_\iota+U_l^i \\
&=&
\sum_{\nu=1}^k\sum_{\alpha=1}^m\partial_\alpha\tilde{u}_\iota^i(X^i)\cdot
e_\alpha^\nu \partial^\gamma
T_\nu^i+\sum_{\alpha=1}^m\partial_\alpha \tilde{u}_\iota^i(X^i)\cdot
S_l^i-\partial^\gamma u_\iota+U_l^i. \eeqnn} \hfill
\parbox{8mm}{\beqn \label{langegleichung} \eeqn} \\
For $1\leq \nu,\iota
\leq k$ we define functions \beqnn A_{\nu,\iota}^i:B_\delta
&\rightarrow& \mathbb{R},\hspace{0.8cm}
A_{\nu,\iota}^{i}:=\sum_{\alpha=1}^{m}\partial_\alpha
\tilde{u}_\iota^i(X^i)\cdot e_\alpha^\nu, \hspace{2cm}\\ \text{and} \hspace{4.55cm} \\
B_{\iota}^{i,\gamma}: B_\delta&\rightarrow& \mathbb{R},
\hspace{0.8cm} B_\iota^{i,\gamma}:=\sum_{\alpha=1}^m\partial_\alpha
\tilde{u}_\iota^i(X^i)\cdot S_l^i-\partial^\gamma u_\iota+U_l^i.
\hspace{2cm} \eeqnn
Inserting these functions into (\ref{langegleichung}) yields
\beqn
\partial^\gamma T_\iota^i=\sum \limits_{\nu=1}^{k}A_{\nu,\iota}^{i}\cdot
\partial^\gamma T_\nu^i+B_{\iota}^{i,\gamma} \;\;\;\;\;\; \text{ for }1\leq \iota
\leq k. \eeqn For each $x \in B_\delta$ we have $\partial^\gamma
T^i(x) \in \mathbb{R}^k$, $B^{i,\gamma}(x) \in \mathbb{R}^k$, and
$A^i(x)=(A_{\nu,\iota}^{i}(x))\in \mathbb{R}^{k\times k}$. We
obtain \beqnn \partial^\gamma T^i=A^i\cdot
\partial^\gamma T^i+ B^{i,\gamma}, \eeqnn hence almost everywhere \beqn \label{opnormestimate}
|\partial^\gamma T^i| \leq
|A^i\cdot
\partial^\gamma T^i|+ |B^{i,\gamma}|\leq \|A^i\|_{\text{op}} |\partial^\gamma T^i|+ |B^{i,\gamma}|,
\eeqn where $\|\cdot\|_{\text{op}}$ denotes the operator norm.
We write $e^\nu=(e^\nu_h,e^\nu_v)\in\R^m\times \R^k$.
By construction of the frame
we have $|e_h^\nu|\leq2\alpha\leq 1$, and by (\ref{alpha1k}) moreover
$\|D\tilde{u}^i\|_{C^{0}(B_{r-\varrho})}\leq \frac{1}{2\sqrt{k}}$. We estimate
\beqnn\|A^i\|_{\text{op}}^2\leq
 \sum_{\nu,\iota=1}^k\left(\sum_{\alpha=1}^m\partial_\alpha\tilde{u}_\iota^i(X^i)\cdot e_\alpha^\nu\right)^{\!2}
=\sum_{\nu,\iota=1}^k\langle D\tilde{u}_\iota^i(X^i),e_h^\nu\rangle^2
\leq \|D\tilde{u}^i(X^i)\|^2\left(\sum_{\nu=1}^k|e_h^\nu|^2\right)\leq \frac{1}{4k}\cdot k=\frac{1}{4}.
\eeqnn Hence with (\ref{opnormestimate}) we obtain \beqn |\partial^\gamma T^i|\leq 2|B^{i,\gamma}|.
\eeqn
By Lemmas \ref{ableitung1} and
\ref{ableitung2}, $B^{i,\gamma}$ only depends on derivatives of
$e^1,\ldots,e^k,u,\tilde{u}^i$ up to the order $|\gamma|$, and on derivatives
of $T^i$ up to the order $|\gamma|-1$. Moreover
$e^1,\ldots,e^k$ are fixed mappings, which are bounded in $C^\mathrm{k}$,
and $u$ is a fixed mapping, which is bounded in
$W^{\mathrm{k},p}$ and hence also in $C^{\mathrm{k}-1}$;
the mappings $\tilde{u}^i$ are uniformly bounded in $W^{\mathrm{k},p}$ and
hence also in $C^{\mathrm{k}-1}$.
From (\ref{vorkoordinate}) it follows that $T^i$ is uniformly bounded in $C^0$. Hence for
$|\gamma|=1$, independently of $i$, we have the estimate
$\|B^{i,\gamma}\|_{L^{\infty}(B_{\delta})}\leq K^1$, and inductively we obtain
for $|\gamma|=\mathpzc{l}$, $1<\mathpzc{l}< \mathrm{k}$,
the estimate $\|B^{i,\gamma}\|_{L^{\infty}(B_{\delta})}\leq
K^{\mathpzc{l}}$, and finally for $|\gamma|=\mathrm{k}$ the estimate
$\|B^{i,\gamma}\|_{L^{p}(B_{\delta})}\leq K^{\mathrm{k}}$. Hence
$T^i$ is uniformly bounded in $W^{\mathrm{k},p}$. As
for $p>m$ the space
$W^{\mathrm{k},p}$ is compactly embedded in $C^{\mathrm{k}-1}$,
there is a function $T:B_\delta \rightarrow \mathbb{R}^k$, and a subsequence of
$T^i$, which converges in
$C^{\mathrm{k}-1}$ to $T$. By Lemma \ref{punktweise} we have $T\equiv 0$. \\
\\ Hence we have locally \beqn f+ \sum
\limits_{\nu=1}^{k}T_\nu^ie^\nu\rightarrow f \text{ in }
C^{\mathrm{k}-1}
\text{ as } i\rightarrow \infty, \eeqn which we wanted to show. \hfill $\square$ \\ \\
Next we like to show that the limit immersion satisfies also the bounds for the
second fundamental form and the volume. For that let $F:\R^N\times \R^{mN}\rightarrow \R$ be a
function, where $N=nm$. For a domain $\Omega\subset \R^m$ and for each $W^{2,1}$-function
$v:\Omega\rightarrow \R^n$ we define
\beqnn
\mathcal{F}(v)=\int_\Omega F(Dv,D^2v)\,d\mathcal{L}^m. \eeqnn
We first need the following lemma:
\begin{lemma} \label{lowsemicont}
Suppose $F:\R^N\times \R^{mN}\rightarrow \R$ is continuous, nonnegative, and
$F(\zeta,\cdot)$ is convex for every fixed $\zeta\in \R^N$. Then, if $v_i,v\in W^{2,1}(\Omega)$ and
$Dv_i\rightarrow Dv$ in $L^1(\Omega)$, $D^2v_i\rightharpoondown D^2v$ weakly in $L^1(\Omega)$, it
follows that \beqnn
\mathcal{F}(v) \leq \liminf_{i\rightarrow \infty}\mathcal{F}(v_i). \eeqnn
\end{lemma}
\textbf{Proof:} \\
This is a special case of Theorem 1.6 in \cite{struwe}. \hfill $\square$ \\ \\
Now we come to the bounds for the limit, which is the final step in the proof of Theorem \ref{compactness0}:
\begin{theorem} The limit immersion $f:M\rightarrow
\mathbb{R}^n$ satisfies $\|A(f)\|_{L^{p}(M)}\leq\mathcal{A}$,
$\vol(M)\leq\mathcal{V}$, and moreover $q \in f(M)$ with $q$ as in Theorem
\ref{compactness0}. \end{theorem}
\textbf{Proof:} \\
We consider the local representations $f\circ \varphi_j^{-1}:B_\delta\rightarrow \R^n$,
$f^i\circ \phi^i\circ \varphi_j^{-1}:B_\delta\rightarrow\R^n$ (where
$\varphi_j:P(B_\delta^j)\rightarrow B_\delta$
is a chart of the atlas $\mathfrak{A}$) and simply write
$f$, $f^i\circ \phi^i$ for that. We consider
the tensorial norm of $A$; for an immersion $f:B_\delta\rightarrow \R^n$
it is pointwisely given by
\beqnn
\|A\|^2=\sum_{i,j,k,l=1}^{m}A_{ij}\cdot A_{kl}\; g^{ik}g^{jl}, \eeqnn
where $G^{-1}=(g^{ij})\in \R^{m\times m}$ is the inverse of $G=Df^t\cdot Df$,
$A_{ij}=(\partial_{ij}f)^\bot$ and $A_{ij}\cdot A_{kl}$ the Euclidean standard scalar product.
Note that the projection onto the normal space only depends on $Df$
(more precisely $\pi^\bot=\Id-DfG^{-1}Df^t$, see for instance \cite{bauer}, p.\ 555).
\\ \\
First we consider the $L^2$-norm of the second fundamental form.
We set \beqnn F=F(Df,D^2f):=\|A\|^2 \eeqnn
and show, that $F$ satisfies the assumptions of Lemma \ref{lowsemicont}. \\ \\
For that we note, that $F$ is a homogeneous polynomial in the variables
$\partial_{ij}f_\mu$ (with $1\leq i,j \leq m$, $1\leq \mu \leq n$)
of degree two, that is a sum of terms of the form $2c\;\partial_{ij}f_\mu\partial_{kl}f_\nu$
and
$c\;(\partial_{ij}f_\mu)^2$, where $c$ only depends on $Df$. \\ \\
Now we write $F=F(\zeta,\xi)$, fix $\zeta$ and calculate the
Hessian $D^2F$ by the variables $\xi$. We obtain
\beqnn D^2F(\zeta,\xi)(v,v)=2F(\zeta,v)\geq 0. \eeqnn
Hence $F$ is convex in $\xi$ for each fixed $\zeta$.
The chain rule implies convexity for the case $p>2$. It remains to show convexity in the
case of dimension $m=1$ and $1<p<2$.
For an immersion $f:(-\delta,\delta)\rightarrow \R^n$,
 the pointwise norm of $A$ simplifies to \beqnn
\|A\|=\frac{1}{|Df|^{2}}|\pi^\bot f''|, \eeqnn
where $\pi^\bot=\Id-DfG^{-1}Df^t$ as above. This time we set $F=F(Df,D^2f):=\|A\|$. Again
we write $F=F(\zeta,\xi)$ and fix $\zeta\neq 0$. Let $c:=\frac{1}{|\zeta|^{2}}$,
$C:=\Id-\frac{1}{|\zeta|^{2}}\zeta\zeta^t$. Then
$F(\zeta,\xi)=c|C\,\xi|$.
Now let $\xi_1,\xi_2 \in \R^n$. Then
for any $t$ with $0\leq t\leq1$, we conclude by the linearity of $C$ that
\beqnn
F(\zeta,(1-t)\xi_1+t\xi_2) &\leq& (1-t)F(\zeta,\xi_1)+tF(\zeta,\xi_2). \eeqnn
Again, this shows that $F$ is convex in $\xi$ for each fixed $\zeta$.
Then it is easily seen, that also \beqnn
\|A\|^p=\frac{1}{|Df|^{2p}}|\pi^\bot f''|^p \eeqnn
is convex in $\xi$ for $p>1$.
In both cases ($p\geq 2$ and $1<p<2$) this implies the desired $L^p$-bound as follows: \\ \\
In the convergence proof it is shown that the local representations of $f^i\circ \phi^i$
are uniformly bounded in $W^{2,p}$. Moreover $f^i\circ \phi^i$ converges in $C^1$ to $f$. As
$W^{2,p}$ is reflexive, there exists a subsequence which converges weakly
in $W^{2,p}$ and therefore also weakly in $W^{2,1}$ to $f$
(see e.g.\ \cite{alt}, p.\ 220, Example 6.10 \nolinebreak 3)\,). Lemma \ref{lowsemicont} gives
\beqnn
\|A(f)\|_{L^p(P(B_\delta^j))}\leq \liminf_{i\rightarrow\infty}
\|A(f^i\circ \phi^i)\|_{L^p(P(B_\delta^j))} \eeqnn
Note that above we have defined $F=F(\zeta,\xi)$ only for $\zeta\in \R^{n\times m}\cong \R^{nm}$ with rank $\zeta=m$;
however Lemma \ref{lowsemicont} is also true under this restriction.
Using a partition of unity we deduce
\beqnn
\|A(f)\|_{L^p(M)}&\leq& \liminf_{i\rightarrow \infty}
\|A(f^i\circ \phi^i)\|_{L^p(M)} \\
&=& \liminf_{i\rightarrow \infty}
\|A(f^i)\|_{L^p(M^i)}\\
&\leq& \mathcal{A}\,. \eeqnn
For the volume we note
$\vol(P(B_\delta^\mathrm{j}))=\int_{B_{\delta}}\sqrt{\det
g_{ij}}\;\;d\mathcal{L}^m$, where $G=(g_{ij})\in \mathbb{R}^{m\times
m}$ with $G=Df^t\cdot Df$.
Now the bound on the volume of the limit manifold $M$ follows directly
from $C^1$-convergence of the local representations.
Finally we note that for each $i$ there is a point $p^i\in M^i$ with
$f^i(p^i)=q$. Then the relation $q \in f(M)$ is obvious. \hfill $\square$ \\
\begin{rem} With the compactness theorem (together with the lower semicontinuity of the norm of the
second fundamental form) it is possible to derive existence theorems for minimizers of the
$L^p$-norm of the second fundamental form, see \cite{langer}, p.\ 224. Analogous results in the setting
of integral rectifiable $m$-varifolds have been attained by J.\ Hutchinson in \cite{hutchinson}
and A.\ Mondino in \cite{mondino}.
\end{rem}
We would like to conclude this section with some generalizations of Theorem \ref{compactness0}. First we would like to show
how to obtain \emph{higher order convergence,} that is convergence in $C^{\K-1}$ for $\K\geq2$\:
(again $\K=$ \!degree of differentiability, $k=$ \!codimension). For that
we assume in addition to (\ref{boundsecond}) the bounds \beqn
\|\nabla^l A(f)\|_{L^\infty(M)}\leq \mathcal{A}_l \hspace{5mm}
\text{ for any } \,l \text{ with }\, 0\,\leq\, l\,\leq\, \K-2, \eeqn
where $\K\geq 2$. Here we assume that each immersion is sufficiently smooth, that is at least of
class $C^{\K}$. Additionally assume that also the mapping $\nu:M\rightarrow G_{n,k}$ we used
to construct the diffeomorphisms $\phi^i$ is at least $C^{\K}$.
For $\alpha>0$ choose $r>0$ such that each immersion is an
$(r,\alpha)$-immersion. We need a bound for higher derivatives of the graph functions $u$:
\begin{lemma} \label{higherderivatives}
For $B_r\subset \R^m$ and $n=m+k$ let $f\in C^\K(B_r,\R^n)$ be a mapping of the type
$f(x)=\nolinebreak(x,u(x))\linebreak\in\R^m \times \R^k$ with $u(0)=0$. Suppose $\|Du\|_{C^0(B_r)}\leq \alpha< \infty$ and
$\|\nabla^l A(f)\|_{L^\infty(B_r)}\leq \mathcal{A}_l<\infty$ for any $l$ with $0\leq l \leq \K-2$. Then \beqnn
\|u\|_{C^\K(B_r)}\,\leq\: C(r,\alpha,\mathcal{A}_0,\ldots,\mathcal{A}_{\K-2}) \eeqnn
for a universal constant $C(r,\alpha,\mathcal{A}_0,\ldots,\mathcal{A}_{\K-2}) < \infty$.
\end{lemma}
\textbf{Proof:} \\
We have $\|Du\|_{C^0(B_r)}\leq \alpha$, hence $\|u\|_{C^0(B_r)}\leq \alpha r$.
The higher derivatives of $u$  are easily estimated by induction, see for instance Lemma 8.2 in the
diploma thesis \cite{schlichting} for
the case of codimension $1$. The proof of the general case is left to the reader.
\hfill $\square$
\\ \\
Starting from the situation in Lemma \ref{higherderivatives}, the calculations in this section show that a subsequence
of $f^i\circ \phi^i$ converges locally in $C^{\K-1}$ to $f$.
We will make use of this when proving higher order convergence in Theorem \ref{compactness1}. \\ \\
Finally we like to explain how to prove the theorem for immersions $f^i:M^i\rightarrow N$ with \emph{values in a
complete Riemannian manifold} $N$ (without boundary). Note that $N$ has to be complete, as otherwise we could
take an open ball $N=B_1\subset\R^n$ and construct a sequence of immersions converging to the boundary $\partial B_1$.
We shall use the Nash embedding: Any Riemannian manifold ($N^n,g$) can be isometrically embedded into $\R^\nu$, where
$\nu=\nu(n)$. Let $\phi:N\rightarrow \R^\nu$ be such an embedding.
Here we additionally assume that the second fundamental form of $\phi$ is bounded in $L^\infty$. \\ \\
Let $f^i:M^i\rightarrow N$ be a sequence of immersions with $\|A(f^i)\|_{L^p(M^i)}\leq \mathcal{A}$,
vol$(M^i)\leq \mathcal{V}$ and $q\in f^i(M^i)$ for a $q\in N$. Now consider the sequence
$\phi\circ f^i:M^i\rightarrow\R^n$. Then $\phi\circ f^i:M^i\rightarrow \R^n$ is a sequence
with $\|A(\phi\circ f^i)\|_{L^p(M^i)}\leq \mathcal{A'}$, vol$(M^i)\leq \mathcal{V}$ and
$\phi(q)\in \phi\circ f^i(M^i)$. Hence we can apply the compactness theorem for immersions into $\R^n$.
We obtain a limit immersion $f:M\rightarrow \R^\nu$ with $f(M)\subset \phi(N)$, and a subsequence of $\phi\circ f^i$
converging to $f$. Applying $\phi^{-1}$ to these mappings, we finally obtain a version of our compactness theorem for immersions
with values in $N$. Again, one can formulate similar statements involving higher order convergence.
\end{section}
\begin{section}{Compactness for immersions on noncompact manifolds}
In the final
section we want to prove Theorem \ref{compactness1}, the
compactness of proper immersions on manifolds
which are not necessarily compact.
\vspace{5mm} \\ \noindent One of the technical main difficulties in the proof
lies in fact that the norm of $A^i$ depends on $R$, that is
\beqnn
\|A^i\|_{L^{\infty}(B_R)} &\leq& C_0(R). \eeqnn
For that reason we do not have uniform estimates for the size of the
radius $r$ of the function graphs as in (\ref{sizegraphradius}) --- we may only estimate $r$ for each fixed $R>0$.
This leads to the problem that we cannot directly apply Lemma \ref{inclusiondipl} b) any more, which was
of great importance for the proof in the compact case.
This explains the need for some technical refinements which are carried out in the following.
\\ \\ \\
\textbf{Preparations for the noncompact case} \\ \\ \\[-4mm]
First of all we have to adjust some definitions to the new situation.
Let $\N$ denote the integers greater than $0$ and let $\N_0=\N\cup\{0\}$.
Any sequence $(a_i)_{i\in\N}$ shall be denoted by $a$.
If $a,b:\N\rightarrow \R$
are two sequences, then $a<b$ if and only if $a_i<b_i$ for all $i\in \N$. \\ \\
Let us assume here that the image $f^i(M^i)$ of any immersion is unbounded in $\R^n$
(as $f^i$ is proper, this is the case whenever $M^i$ is noncompact). If the $f^i(M^i)$ are bounded
uniformly in $i$,
then the proof of Theorem \ref{compactness0} applies. If the $f^i(M^i)$ are bounded, but
not uniformly in $i$, then the statement is proven as in the unbounded case but with minor adaptations
of the notation (for a more general formulation see Corollary \ref{coropen}). \\ \\
We are dealing here with balls both in $\R^m$ and in $\R^n$. For $\varrho>0$ let $B_\varrho$
denote the open ball of radius $\varrho$ in $\R^m$ centered at the origin.
Let $\hat{B}_\varrho$ denote the corresponding ball
in $\R^n$. For $\varrho\leq 0$ we define $\hat{B}_\varrho=\emptyset$.
Note that all balls $B_R$ and $B_i$ in Theorem \ref{compactness1}
are in fact balls in $\R^n$ and should be written
as $\hat{B}_R$ and $\hat{B}_i$ in our new notation. \\ \\
For a given immersion $f:M\rightarrow \R^n$ and for $p\in M$ we set \beqnn
\bar{p}=\bar{p}\,(f):=\min\{j\in \N:f(p)\in \hat{B}_j\} \in \N. \eeqnn
The notion of the $(r,\alpha)$-immersion has to be adapted by replacing the real number $r$ by
a sequence. First we adapt the definition of $U_{r,q}$:
\begin{defi} \label{decomp1} Let $f:M\rightarrow \R^n$ be an immersion and let $q\in M$. Let $A_q$ and $\pi$ be
as in Section \ref{locrepfg}. Let $r:\N\rightarrow \R_{>0}$ be a sequence. We define $U_{r,q}$ to be the
$q$-component of the set $(\pi\circ A_q^{-1}\circ f)^{-1}(B_{r_{\bar{q}}})$. \end{defi}
With the preceding definition we come to the notion of an $(r,\alpha)$-immersion:
\begin{defi} \label{ralphanoncompact} Let $f:M\rightarrow \R^n$ be an immersion. Let $r:\N\rightarrow \R_{>0}$
be a sequence and let $\alpha>0$.
We say that $f$ is an $(r,\alpha)$-immersion, if for each
$q\in M$ the set $A_q^{-1}\circ f(U_{r,q})$ is the graph of
a $C^1$-function $u:B_{r_{\bar{q}}}\rightarrow \R^k$ with $\|Du\|_{C^{0}(B_{r_{\bar{q}}})}\leq \alpha$.
\end{defi}
Under the assumption that each immersion is proper,
the condition
$\|A^i\|_{L^{\infty}(\hat{B}_R)} \leq C_0(R)$ obviously implies, that
for every $\alpha>0$ there
is a sequence $r$ (which does not depend on $i$) such that each immersion $f^i$ is an
$(r,\alpha)$-immersion. From now on $r$ will always be a sequence. All sequences
$r$ and $\varrho$ with $\varrho\leq r$ are assumed to be greater than $0$.
\begin{defi} \label{decomp2} $ $ \\ \vspace{-5.3mm}
\begin{itemize}
\item[a)] Let $\nu:\N_0\rightarrow \N_0$ be a sequence.
We say $\nu$ is a subdivision, if $\nu_0=0$ and
$\nu$ is strictly increasing.
\item[b)] Let $f:M\rightarrow \R^n$ be
an $(r,\alpha)$-immersion, $\nu$ a subdivision and $\delta$
a sequence with $\delta<r$. Let
$Q=\{q_1,q_2,\ldots \}$ be a countable set of points in $M$. We say
$Q$ is a $\delta$-net for $f$ with subdivision \nolinebreak$\nu$, if for all $j
\in \N$ the following holds:
\begin{itemize}
\item[$\bullet$] $f(q_k)\in \hat{B}_j\setminus \hat{B}_{j-1}$ for all $k$ with
$\nu_{j-1} < k \leq \nu_j$,
\item[$\bullet$] $f^{-1}(\hat{B}_j)\subset \bigcup_{k=1}^{\nu_{j}}U_{\delta,q_{k}}$.
\end{itemize} \end{itemize} \end{defi} \vspace{3mm}
Next we have to adapt the definition of $(\mathfrak{G}^s,\mathfrak{d})$
in order to handle graphs with different radii:
\begin{defi} Let $r:\N\rightarrow \R_{>0}$ be a decreasing sequence, $\nu$ a subdivision
and $\varrho:\N\rightarrow \R_{>0}$ a sequence with $\varrho_i=r_k$ for all $i,k\in \N$
with $\nu_{k-1}<i\leq \nu_k$.
For $j\in \N\cup \{\infty\}$ with $\nu_\infty:=\infty$ we set \\
\parbox{12cm}{\beqnn \mathfrak{G}^j=\mathfrak{G}^j(r,\nu)=\{(A_i,u_i)_{i=1}^{\nu_{j}}: &A_i&:
\mathbb{R}^n\rightarrow \mathbb{R}^n \text{ is a Euclidean isometry,} \\ &u_i& \in
C^1(\overline{B}_{\varrho_i},\mathbb{R}^k)\}. \eeqnn }
\hfill\parbox{8mm}{\beqn \label{grund} \eeqn} \\
If \:$\tilde{\Gamma}=(A_i,u_i)_{i=1}^\infty \in \mathfrak{G}^\infty$, we define
$\tilde{\Gamma}_l:=(A_i,u_i)_{i=1}^{\nu_{l}}\in \mathfrak{G}^l$.
Moreover, splitting $A_i$ in a rotation $R_i\in \mathbb{SO}(n)$ and a translation
$T_i\in \R^n$, we set for $j\in \N$
\\ \parbox{12cm}{ \beqnn \mathfrak{d}(\cdot,\cdot): \mathfrak{G}^j\times
\mathfrak{G}^j&\rightarrow& \mathbb{R}, \\
\mathfrak{d}(\Gamma,\tilde{\Gamma})&=&\sum
\limits_{i=1}^{\nu_j}(\|R_i-\tilde{R}_i\|+|T_i-\tilde{T}_i|+
\|u_i-\tilde{u}_i\|_{C^{1}({B}_{\varrho_{i}})}).
\eeqnn} \hfill \parbox{8mm}{\beqn \label{metrikgraph} \eeqn} \end{defi}
More accurately we should write $\mathfrak{d}^j$ instead of \,$\mathfrak{d}$, but we will
maintain the notation without index $j$.
Again, as in the compact case, $(\mathfrak{G}^j,\mathfrak{d})$ is a metric space.
The reader should take care not to confuse the isometries $A_j$ (of a fixed immersion)
with the second fundamental
forms $A^i$ of the sequence of immersions $f^i$.
\begin{defi} Let $f:M\rightarrow \R^n$ be an $(r,\alpha)$-immersion,
$\delta<r$ a sequence and
$\nu$ a subdivision. Let $Q=\{q_1,q_2,\ldots\}$ be a
$\delta$-net for $f$ with subdivision $\nu$.
As in the compact case we may assign to each $q_j\in Q$
a neighborhood $U_{r,q_{j}}$, a Euclidean isometry $A_j$
and a $C^\infty$-function $u_j:B_{r_{\bar{q}_{j}}}\rightarrow \R^k$.
We define \beqn \Gamma=\Gamma(f)=(A_j,u_j)_{j=1}^{\infty} \in
\mathfrak{G}^\infty. \label{graphsystinfinite} \eeqn Furthermore for each $j \in \N$ we
define \beqn Z(j):=\{l\in\N: U_{\delta,q_{j}}\cap
U_{\delta,q_{l}}\neq \emptyset\}\,\in \mathcal{P}(\N). \label{intersectioninfinite} \eeqn
\end{defi}
\vspace{3mm}
\begin{lemma} \label{distancedeltanet} Let $f:M\rightarrow \R^n$ be an $(r,\alpha)$-immersion
with $\alpha^2 < \frac{1}{3}$ and $r_i\leq \frac{3}{8}$ for all $i\in \N$.
Let $\delta$ be a sequence with $\delta<r$ and let $j \in \N$.
\begin{itemize}
\item[a)] If $p\in f^{-1}(\hat{B}_j\setminus \hat{B}_{j-1})$, then
$U_{\delta,p}\subset f^{-1}(\hat{B}_{j+\frac{1}{2}}\setminus \hat{B}_{j-\frac{3}{2}})$.
\item[b)]
If $p \in f^{-1}(\hat{B}_j)$ and $q \in
M\setminus f^{-1}(\hat{B}_{j+1})$, then $U_{\delta,p}\cap
U_{\delta,q}=\emptyset$.
\end{itemize} \end{lemma} \noindent \textbf{Proof:}
\begin{itemize}
\item[a)]
Let $x \in U_{\delta,p}$. With Lemma \ref{inclusiondipl} a)
we calculate \beqnn
|f(x)-f(p)|\;\leq\; (1+\alpha^2)\,\delta_j\;\leq\; \frac{1}{2}. \eeqnn
As $j-1\leq |f(p)|<j$ it follows $j-\frac{3}{2}\leq |f(x)|< j+\frac{1}{2}$. Hence
$U_{\delta,p}\subset f^{-1}(\hat{B}_{j+\frac{1}{2}}\setminus \hat{B}_{j-\frac{3}{2}})$.
\item[b)] By part a) we have $U_{\delta,p}\subset f^{-1}(\hat{B}_{j+\frac{1}{2}})$,
$U_{\delta,q}\subset f^{-1}(\R^n\setminus \hat{B}_{j+\frac{1}{2}})$, hence
$U_{\delta,p}\cap U_{\delta,q}=\emptyset$. \hfill $\square$
\end{itemize}
\vspace{3mm}
We obtain the following version of Lemma \ref{inclusiondipl} b):
\begin{lemma} \label{inclusionnoncompact} Let $f:M\rightarrow \R^n$ be an $(r,\alpha)$-immersion with
$\alpha^2< \frac{1}{3}$, where $r$ is a decreasing sequence with
$r_1\leq \frac{3}{8}$.
Let $\delta'$ be a decreasing sequence with
$\delta'_i\leq \frac{r_{i+1}}{4}$ for all
$i\in\N$ and let $\delta$ be a sequence with
$\delta_i\leq \frac{\delta'_{i+1}}{4}$ for all $i\in \N$.
\begin{itemize}
\item[a)] If $p,q\in M$ and $U_{\delta,q}\cap U_{\delta,p}\neq
\emptyset$, then
$U_{\delta,p}\subset U_{\delta',q}\subset U_{r,q}$.
\item[b)] If $p,q\in M$ and $U_{\delta',q}\cap U_{\delta',p}\neq
\emptyset$, then
$U_{\delta',p}\subset U_{r,q}$.
\item[c)] If $x,y,z \in M$ and $U_{\delta,x}\cap U_{\delta,y} \neq
\emptyset$, $U_{\delta,y} \cap U_{\delta,z}\neq \emptyset$, then
$U_{\delta,z}\subset U_{r,x}$.
\end{itemize}
\end{lemma} \pagebreak
\textbf{Proof:} \begin{itemize}
\item[a)] Let $j:=\bar{p}$, $k:=\bar{q}$, in other words $p \in f^{-1}(\hat{B}_j\setminus \hat{B}_{j-1})$,\,
$q\in f^{-1}(\hat{B}_k\setminus \hat{B}_{k-1})$. As $U_{\delta,q}\cap U_{\delta,p}\neq
\emptyset$, Lemma \ref{distancedeltanet} b) implies $|j-k|\leq 1$.
Let $\iota=\min\{j,k\}$ and let $\zeta\in U_{\delta,p}$,\, $\xi \in U_{\delta,q}\cap U_{\delta,p}$.
With $\varphi_q=\pi\circ A_q^{-1}\circ f$, using $\delta_i\leq \frac{\delta'_{i+1}}{4}$
and the fact that $\delta'$ is decreasing, we estimate
\beqnn
|\varphi_q(\zeta)|&\leq& |f(\zeta)-f(q)| \\
&\leq& |f(\zeta)-f(p)|+|f(p)-f(\xi)|+|f(\xi)-f(q)| \\
&\leq& 3(1+\alpha^2)\frac{\delta'_{\iota+1}}{4} \\
&<& \delta'_{\iota+1}.
\eeqnn
Hence $U_{\delta,p}\subset \varphi_q^{-1}(B_{\delta'_{\iota+1}})$.
As $\bar{q}\in \{\iota\,,\,\iota+1\}$ and as $\delta'$ is decreasing, we
conclude $U_{\delta,p}\subset \varphi_q^{-1}(B_{\delta'_{\bar{q}}})$. But
$U_{\delta,p}\cup U_{\delta,q}$ is a connected set containing $q$ and is hence
included in the $q$-component of the set $\varphi_q^{-1}(B_{\delta'_{\bar{q}}})$,
that is in $U_{\delta',q}$. Hence $U_{\delta,p}\subset U_{\delta',q}$. The
relation $U_{\delta',q}\subset U_{r,q}$ is obvious.
\item[b)] The proof of the second part runs as before.
\item[c)] As $\delta<\delta'$, the relation $U_{\delta,x}\cap U_{\delta,y} \neq
\emptyset$ implies $U_{\delta',x}\cap U_{\delta',y} \neq \emptyset$. By part a) we have
$U_{\delta,z}\subset U_{\delta',y}$ and by part b) $U_{\delta',y}\subset U_{r,x}$.
\hfill $\square$
\end{itemize} \vspace{0.55mm}
\begin{rem} Let $f:M\rightarrow \R^n$ be an $(r,\alpha)$-immersion with $\alpha^2<\frac{1}{3}$,
$\delta<r$ a sequence and $p,q \in M$ with
$U_{\frac{\delta}{4},q}\cap U_{\frac{\delta}{4},p}\neq \emptyset$. Then, under the
additional assumption $p,q \in f^{-1}(\hat{B}_j\setminus \hat{B}_{j-1})$
for a $j\in \N$, we may apply Lemma \ref{inclusiondipl} b) and obtain
\beqn \label{inclusionspecial} U_{\frac{\delta}{4},p}\subset U_{\delta,q}. \eeqn
This will be used in the proof of the following lemma.
\end{rem} \vspace{0.2mm}
\begin{lemma} \label{existencedeltanet}
Let $f^i:M^i\rightarrow \R^n$ be a sequence as in Theorem \ref{compactness1}.
Moreover let $\alpha>0$ with $\alpha^2< \frac{1}{3}$ and
$r$ a sequence with $r_i\leq \frac{3}{8}$ for all $i\in \N$, such that each $f^i$ is an
$(r,\alpha)$-immersion. Let $\delta$ be a sequence with
$\delta<r$.
Then there exists
a fixed subdivision $\nu$, such that the following holds:
\begin{itemize}
\item[a)] Each immersion
$f^i$ admits a $\delta$-net $Q^i$ with
subdivision $\nu$.
\item[b)] For
each $j \in \N$ let $Z^i(j)$ be the set corresponding to $Q^i$ as defined
in (\ref{intersectioninfinite}). Then, after passing to a subsequence,
for each $j \in \N$ there exists a finite set
$Z(j)\subset \N$, such that
\beqnn Z^i(k)=Z(k) \text{ for all } k\leq \nu_i. \eeqnn
\end{itemize}
\end{lemma}
\noindent \textbf{Proof:}
\begin{itemize}
\item[a)] We fix $i,j \in \N$. Now consider the immersion $f^i$. Let
$q_1^j \in (f^i)^{-1}(\hat{B}_j \setminus \hat{B}_{j-1})$. Assume we
have found points $\{q_1^j, \ldots, q_\iota^j\}$ in $(f^i)^{-1}(\hat{B}_j
\setminus \hat{B}_{j-1})$ with the property $U^i_{\delta/4,q_a^j}\cap
U^i_{\delta/4,q_b^j}=\emptyset$ for $a \neq b$. Suppose
$U^i_{\delta,q_1^j}\cup \ldots \cup U^i_{\delta,q_\iota^j}$ does not cover
$(f^i)^{-1}(\hat{B}_j\setminus \hat{B}_{j-1})$. Then choose a point
$q_{\iota+1}^j\in (f^i)^{-1}(\hat{B}_j\setminus \hat{B}_{j-1})$ from the complement. Then
$U^i_{\delta/4,q_k^j}\cap U^i_{\delta/4,q_{\iota+1}^j}=\emptyset$ for $k\leq \iota$, as
otherwise $U^i_{\delta/4,q_{\iota+1}^j}\subset U^i_{\delta,q_k^j}$ by (\ref{inclusionspecial}).
Using (\ref{comp1eq1}) in the first line and
Lemma \ref{distancedeltanet} a) in the second, we estimate \beqnn
C(j+1) &\geq& \mu^i(\hat{B}_{j+\frac{1}{2}}\setminus \hat{B}_{j-\frac{3}{2}}) \\
&\geq& \sum_{k=1}^s \mu_{g^i}(U^i_{\delta/4,q_k^j}) \\
&\geq& \sum_{k=1}^s \mathcal{L}^m(B_{\frac{\delta_j}{4}}) \\
&\geq& s\left(\frac{\delta_j}{4}\right)^m. \eeqnn
Therefore, with $\lfloor x \rfloor:=\max\{n\in\N_0:n\leq x\}$ for $x\geq 0$,
this procedure yields after at most $\lfloor(\frac{4}{\delta_j})^mC(j+1)\rfloor$ steps a cover
of $(f^i)^{-1}(\hat{B}_j\setminus \hat{B}_{j-1})$. Now define the subdivision
$\nu$ recursively as follows: \beqnn
\nu_0:=0, \hspace{1.3cm} \nu_j:=\nu_{j-1}+\Bigl\lfloor\Bigl(\frac{4}{\delta_j}\Bigr)^mC(j+1)\Bigr\rfloor \hspace{5mm}\text{ for }\,
j\geq 1. \eeqnn
By the considerations of above we may choose for all $i,j\in \N$
exactly $\nu_j-\nu_{j-1}$ points
$q_{\nu_{j-1}+1}^i,\ldots,q_{\nu_j}^i$ in $(f^i)^{-1}(\hat{B}_j\setminus \hat{B}_{j-1})$,
such that $(f^i)^{-1}(\hat{B}_j)\subset \bigcup_{k=1}^{\nu_j}U^i_{\delta,q_k^i}$.
\item[b)] Fix $j \in \N$. Let $k \leq \nu_j$. If $l>\nu_{j+1}$ then
$U^i_{\delta,q_k^i}\cap U^i_{\delta,q_l^i}=\emptyset$ by Lemma \ref{distancedeltanet} b), hence $l
\notin Z^i(k)$ for each $i\in \N$. \!But this means
$Z^i(k)\!\subset\! \{1,\ldots, \nu_{j+1}\}$, hence $|\{Z^i(k):i\in\nolinebreak
\N\}|\leq |\mathcal{P}(\{1, \ldots,
\nu_{j+1}\})|=2^{\nu_{j+1}}$. Hence we may pass to a subsequence
$(f^{a_{i}^{j}})_{i\in\N}$ of $(f^i)_{i\in\N}$ with \beqnn Z^{a_{i}^{j}}(k)=Z(k) \text{
for all } k\leq \nu_j \text{ and all } i \in \N. \eeqnn
Choosing successively subsequences for any $j\in \N$ and passing to the diagonal sequence,
we obtain a subsequence with the desired property. (Note that the sets $Z(k)$ do not depend
on $j$,\, if passing successively to subsequences.)
 \hfill $\square$
\end{itemize}
\vspace{3mm}
\begin{lemma} \label{convergencegraphsystem} Let $f^i:M^i\rightarrow \R^n$ be a sequence as in Theorem
\ref{compactness1} and $r$ a sequence, such that each immersion $f^i$ is an $(r,\alpha)$-immersion.
Let $\delta<r$ be another sequence and $\nu$ a subdivision.
Let $Q^i$ be $\delta$-nets for $f^i$ with subdivision $\nu$ and let \,$\Gamma^i\in \mathfrak{G}^\infty$
be as in (\ref{graphsystinfinite}).
Then, after passing to a subsequence, there exists a graph system $\Gamma \in
\mathfrak{G}^{\infty}$, such that
for all $j\in \N$ \beqnn \Gamma_j^i\rightarrow
\Gamma_j \;\;\text{ in } (\mathfrak{G}^{j},\mathfrak{d})\;
\text{ as } i\rightarrow\infty. \eeqnn
\end{lemma}
\noindent \textbf{Proof:} \\
Fix $j \in \N$. With the arguments of Theorem 3.3 in \cite{langer},
there exists a graph system $\tilde{\Gamma}^j\in
\mathfrak{G}^j$ and a subsequence $(f^{a_{i}^{j}})_{i\in\N}$
of $(f^i)_{i\in\N}$, such that \beqnn
\Gamma_j^{a_{i}^{j}}\rightarrow \tilde{\Gamma}^j \text{ in }
(\mathfrak{G}^{j},\mathfrak{d}) \text{ as }
i\rightarrow\infty. \eeqnn
By successively choosing subsequences for any $j\in \N$ and passing to the diagonal sequence,
we obtain a sequence with $\Gamma_j^i\rightarrow \tilde{\Gamma}^j$ in
$(\mathfrak{G}^{j},\mathfrak{d})$ as $i\rightarrow\infty$ for all $j\in\N$.
Moreover, if $\tilde{\Gamma}^k=(\hat{A}_j^k,\hat{u}_j^k)_{j=1}^{\nu_{k}}$ and
$\tilde{\Gamma}^l=(\hat{A}_j^l,\hat{u}_j^l)_{j=1}^{\nu_{l}}$ with
$k\leq l$, we observe
$(\hat{A}_j^k,\hat{u}_j^k)=(\hat{A}_j^l,\hat{u}_j^l)$ for $j \leq
\nu_k$. Define $\Gamma=(A_j,u_j)_{j=1}^\infty \in
\mathfrak{G}^\infty$ by setting
$(A_j,u_j):=(\hat{A}_j^k,\hat{u}_j^k)$ for an arbitrary $k$ with
$\nu_k\geq j$. Hence $\Gamma_j=\tilde{\Gamma}^j$ for each $j \in
\N$, which completes the proof. \hfill $\square$
\\ \\ \\
\textbf{Construction of the limit manifold and immersion} \\ \\ \\[-4mm]
Let $f^i:M^i\rightarrow \R^n$ be a sequence of
immersions as in Theorem \ref{compactness1}.
All constants have to be chosen such that all arguments of the compact case can be used.
Let $\alpha>0$ with
$\alpha\leq \frac{1}{4\sqrt{k}}$ as in (\ref{alpha1k}),
in particular $\alpha^2\leq\frac{1}{10}$ as in Section \ref{secrepimm}.
Let $r$ be a decreasing sequence with $r_1\leq \frac{3}{8}$, such that each $f^i$ is an
$(r,\alpha)$-immersion. Let the sequence $\delta'$ be defined by
$\delta'_i=\frac{r_{i+1}}{4}$ for all $i \in \N$, and $\delta$ be defined by
$\delta_i=\frac{r_{i+2}}{16}$ for
all $i\in \N$, that is $\delta_i=\frac{\delta'_{i+1}}{4}$.
Finally let $\hat{\delta}$ be the sequence defined by
$\hat{\delta}_i=\delta_{i+1}$ for all $i\in \N$.
By Lemma \ref{existencedeltanet} a) there
exist a fixed subdivision $\nu$ and countable subsets
$Q^i=\{q_1^i,q_2^i,\ldots \} \subset M^i$,
such that each $Q^i$ is a $\frac{\hat{\delta}}{10}$-net for $f^i$ with
subdivision $\nu$. Similar to the compact case, we have to use $\frac{\hat{\delta}}{10}$-nets
and not only $\delta$-nets (in the compact case we used $\frac{\delta}{10}$-nets).
Moreover, by Lemma \ref{existencedeltanet} b), we may pass to a subsequence such that for each $j \in \N$ there
exists a finite set $Z(j)\subset \N$
with \beqn \label{ZikZk2} Z^i(k)=Z(k) \text{ for all } k \leq \nu_i.
\eeqn
(Here $Z^i(k):=\{l\in\N: U_{\delta,q_{k}^i}^i\cap
U_{\delta,q_{l}^i}^i\neq \emptyset\}$ as in (\ref{intersectioninfinite}),
that is $Z^i(k)$ is
\emph{not} defined as $\{l\in\N: U_{\hat{\delta}/10,q_{k}^i}^i\cap
U_{\hat{\delta}/10,q_{l}^i}^i\neq \emptyset\}$. Nevertheless, as any
$\frac{\hat{\delta}}{10}$-net is also a $\delta$-net, (\ref{ZikZk2}) holds.) \\ \\
By Lemma \ref{convergencegraphsystem}, after passing to another subsequence,
there exists a graph system $\Gamma=(A_i,u_i)_{i=1}^\infty \in
\mathfrak{G}^{\infty}$, such that
for each $j \in \N$ \beqnn \Gamma_j^i\rightarrow
\Gamma_j \text{ in } (\mathfrak{G}^{j},\mathfrak{d})
\text{ as } i\rightarrow\infty. \eeqnn
\\ We come to the construction of the limit manifold and limit immersion: \\ \\
Let $\rho$ be a sequence with $\rho_j=\delta_k$ for all $j,k\in \N$ with
$\nu_{k-1}<j\leq \nu_k$. We define
$B_\delta^j:=B_{\rho_{j}}\times \{j\}$.
The set $\bigcup_{j=1}^\infty B_\delta^j$, endowed with the
disjoint union topology, is a second countable space. Again we define
a relation $\sim$ on $\bigcup_{j=1}^\infty B_\delta^j$. For
$(x,j),(y,k)\in \bigcup_{l=1}^\infty B_\delta^l$ we set \beqnn (x,j)
\sim (y,k)\; \Leftrightarrow\; [k \in Z(j) \text{ and }
A_j(x,u_j(x))=A_k(y,u_k(y))].  \eeqnn
Here the sets $Z(j)$ shall be the fixed sets from (\ref{ZikZk2}). \\ \\
To simplify the notation, for any sequence $\varrho$ with $0<\varrho\leq r$ we set \beqnn
U_{\varrho,j}^i:=U_{\varrho,q_{j}^i}^i. \eeqnn
Now observe that the construction of the limit manifold $M$ can be
performed in exactly the same manner as in the compact case.
For that we note that the sequence $r$ in the present case corresponds to the number
$r$ in the compact case. Similarly, the sequence $\delta'$ corresponds to the number $\frac{r}{4}$,
the sequence $\delta$ to the number $\delta=\frac{r}{16}$, the sequence
$\frac{\hat{\delta}}{10}$ to the number
$\frac{\delta}{10}$. \\ \\
Lemma \ref{inclusiondipl} is replaced by Lemma
\ref{inclusionnoncompact}. By part c) of the latter, even the iterated case works.
Moreover, in the compact case it was crucial that the sets $Z^i(k)$ do not depend on $i$,
that is $Z^i(k)=Z(k)$. This is replaced by Lemma \ref{existencedeltanet} b), which
ensures for fixed $k\in \N$ that $Z^i(k)=Z(k)$ for $i$ sufficiently large. In the compact case, all arguments
involving $Z^i(k)=Z(k)$ were either needed for the construction of the limit,
for which it is sufficient to consider $i$ large, or for the reparametrizations
$\phi^i:M\rightarrow M^i$, which are in the present case replaced by diffeomorphisms
$\phi^i:U^i\rightarrow (f^i)^{-1}(B_i)\subset M^i$ for which the property
(\ref{ZikZk2}) suffices. For the same reasons the convergence of graph systems
$\Gamma_j^i\rightarrow \Gamma_j$ for any $j\in \N$,
replacing $\Gamma^i\rightarrow \Gamma$, is sufficient for our proof. \\ \\
Following step by step the arguments of Lemma \ref{equivalence}, we see
that $\sim$ defines an equivalence relation on $\bigcup_{j=1}^\infty B_\delta^j$. Again we set
$M=(\bigcup_{j=1}^\infty B_\delta^j)/\!\sim$\,. We construct an atlas $\mathfrak{A}$ as
in Section \ref{constructionsection}, with charts $\varphi_{_V}^j:P(V^j)\rightarrow V$ for
$V\subset B_{\rho_j}$. \\ \\
Similarly, we may follow the arguments of Lemmas \ref{quotientopen}, \ref{MisHausdorff},
and \ref{diffatlas},
stating that the quotient projection $P:\bigcup_{j=1}^\infty B_\delta^j\rightarrow M$ is open,
that $M$ is a second countable Hausdorff space and $\mathfrak{A}$ a differentiable atlas on
$M$. Hence $(M,\mathfrak{A})$ induces uniquely the structure of a differentiable manifold. \\ \\
Finally we define a smooth immersion on $M$ by
\beqnn f:\hspace{6mm}M\hspace{3.5mm}&\rightarrow& \R^n, \\
\text{[}(x,j)\text{]}&\mapsto& A_j(x,u_j(x)), \eeqnn where $[(x,j)]$ denotes the
equivalence class of $(x,j)$. \pagebreak \\
We have the following versions
of Lemmas \ref{eigenschaft1} and \ref{eigenschaft2}:
\begin{itemize}
\item It holds $M=\bigcup_{j=1}^\infty P(B_{\hat{\delta}/6}^j)$.
\item If\, $P(B_{\hat{\delta}/4}^j)\cap P(B_{\hat{\delta}/4}^k)\neq \emptyset$, then
$P(B_{\hat{\delta}/4}^k)\subset P(B_\delta^j)$.
\end{itemize} \vspace{3mm}
As in the compact case, let us define sets \beqnn
\tilde{Z}^i(j)=\{l\in \N:\, U_{\hat{\delta}/5,j}^i\cap U_{\hat{\delta}/5,l}^i\neq \emptyset\}.
\eeqnn
With the arguments of Lemma \ref{existencedeltanet} b), we may pass to a subsequence, such that
\beqnn
\tilde{Z}^i(k)=\tilde{Z}(k) \text{ for all } k\leq \nu_i \eeqnn
for fixed finite sets $\tilde{Z}(k)\subset \N$. \\ \\
Finally, the following version of Lemma \ref{eigenschaft3} holds: \begin{itemize}
\item If\, $P(B_{\hat{\delta}/4}^j)\cap P(B_{\hat{\delta}/4}^k)= \emptyset$, then
$k \notin \tilde{Z}(j)$.
\end{itemize} \vspace{2mm}
As additional lemma we have \begin{lemma} The immersion
$f:M\rightarrow \R^n$ is proper. \end{lemma} \noindent
\textbf{Proof:} \\
Let $K\subset \R^n$ be compact. Then there is a $j\in \N$ with $K\subset \hat{B}_j$.
Let $x\in f^{-1}(\hat{B}_j)$. As $M=\bigcup_{l=1}^\infty P(B_{\delta/2}^l)$, there is a
$k\in \N$ with $x\in P(B_{\delta/2}^k)$. It holds
$f(P(B_{\delta/2}^k))=A_k(\{(y,u_k(y)):\nolinebreak y\in \nolinebreak B_{\rho_k/2}\})$. With the argument
of Lemma \ref{distancedeltanet} a)
we conclude $f(P(B_{\delta/2}^k))\subset \hat{B}_{j+1}$. As
$f^i(q_k^i)=A_k^i(0,u_k^i(0))\rightarrow A_k(0,u_k(0))$, we conclude $f^i(q_k^i)\in \hat{B}_{j+1}$
for $i$ sufficiently large. As $Q^i$ is a $\delta$-net for $f^i$ with subdivision $\nu$, we conclude
$k\leq \nu_{j+1}$. Hence \beqnn f^{-1}(K)\subset
f^{-1}(\hat{B}_j)\subset \bigcup_{l=1}^{\nu_{j+1}}P(B_{\delta/2}^l)
\subset \bigcup_{l=1}^{\nu_{j+1}}P(\overline{B_{\delta/2}^l}), \eeqnn
that is $f^{-1}(K)$ is a subset of a compact set. As $f$ is an immersion, it is also continuous,
hence $f^{-1}(K)$ closed. But closed subsets of compact sets are compact.
\hfill $\square$ \\ \\
Note that the limit manifold $M$ does not need to be connected, even
if all manifolds $M^i$ are connected. A simple counterexample is given
in Figure \ref{limitnotconnected}: \\[-5mm]
\setlength{\unitlength}{1cm}
\begin{picture}(10,3.5)
\put(0,0){\line(1,0){2}}
\put(0,2){\line(1,0){2}}
\qbezier(2,2)(3.5,2)(3.5,1)
\qbezier(2,0)(3.5,0)(3.5,1)
\put(3.5,2){\line(1,0){2}}
\put(3.5,0){\line(1,0){2}}
\qbezier[40](2,2)(2.75,2)(3.5,2)
\qbezier[40](2,0)(2.75,0)(3.5,0)
\qbezier(5.5,2)(7,2)(7,1)
\qbezier(5.5,0)(7,0)(7,1)
\put(7,2){\line(1,0){2}}
\put(7,0){\line(1,0){2}}
\qbezier[40](5.5,2)(6.25,2)(7,2)
\qbezier[40](5.5,0)(6.25,0)(7,0)
\qbezier(9,2)(10.5,2)(10.5,1)
\qbezier(9,0)(10.5,0)(10.5,1)
\put(10.5,2){\line(1,0){2}}
\put(10.5,0){\line(1,0){2}}
\qbezier[40](9,2)(9.75,2)(10.5,2)
\qbezier[40](9,0)(9.75,0)(10.5,0)
\qbezier(12.5,2)(14,2)(14,1)
\qbezier(12.5,0)(14,0)(14,1)
\put(2,1){$f^1$}
\put(2.5,1.2){\vector(2,1){0.65}}
\put(5.5,1){$f^2$}
\put(6,1.2){\vector(2,1){0.65}}
\put(9,1){$f^3$}
\put(9.5,1.2){\vector(2,1){0.65}}
\put(12.5,1){$f^4$}
\put(13,1.2){\vector(2,1){0.65}}
\put(14.2,1){$\ldots \,f^i \ldots$}


\end{picture} \\[-4mm]
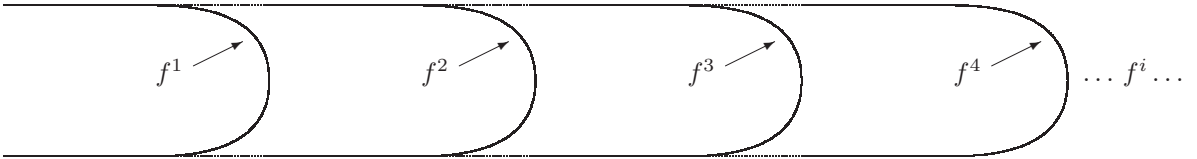
\begin{fig} \label{limitnotconnected} The limit manifold $M$
does not need to be connected, even if all $M^i$ are connected. The limit $f(M)$ in
the example are two parallel lines. \end{fig}
\vspace{5mm}
Let $\nu_f:M\rightarrow G_{n,k}$ denote the Gauss normal map with respect to
the immersion $f$. With Lemma \ref{higherderivatives} we may conclude (as in Step 6 below) that
the limit $f$ is in $C^\infty$. Hence also $\nu_f$ is in $C^\infty$.
We come to the proof of our second main theorem.
\\ \\ \\[-2mm]
\noindent \textbf{Proof of Theorem \ref{compactness1}:} \\ \\[2mm]
\noindent \textbf{Step 1:\; Definition of maps \boldmath $\varphi^i$ \unboldmath} \\ \\
First we define maps
$\varphi^i:\bigcup_{j=1}^{\nu_{i}}P(B_\delta^j)\rightarrow
\bigcup_{j=1}^{\nu_{i}} U_{r,j}^i$ as follows: \\ \\
Let $i \in \N$ be fixed. With the arguments of the compact case,
we may choose $a_i\in\N$ sufficiently large with $a_i\geq i$, such that
for all $j\leq \nu_i$ and all $x \in
P(B_\delta^j)$ the affine space $h(x):=f(x)+\nu_f(x)$ intersects
$f^{a_{i}}(U_{r,j}^{a_{i}})$ in exactly one point $S_x$ and that
this point lies in $f^{a_{i}}(U_{\delta',j}^{a_{i}})$.
Furthermore there is exactly one point $\sigma_x \in
U_{r,j}^{a_{i}}$ with $f^i(\sigma_x)=S_x$. We define
\beqnn \varphi^{a_i}:\;\;\:\bigcup_{j=1}^{\nu_{i}}P(B_\delta^j)&\rightarrow&
\bigcup_{j=1}^{\nu_{i}} U_{r,j}^{a_{i}}\;\,\subset\: M^{a_{i}}\!, \\[1mm]
x\hspace{3mm}&\mapsto&\sigma_x.
\eeqnn \\[-3mm]
We like to show that $\varphi^{a_i}$ is well-defined: \:
Suppose $x$ also lies in $P(B_\delta^k)$ for a $k\leq \nu_i$.
Then
$h(x)$ intersects $f^{a_{i}}(U_{r,k}^{a_{i}})$ in exactly
one point $S_x'$ and there is exactly one $\sigma_x'\in U_{r,k}^{a_{i}}$
with $f(\sigma_x')=S_x'$. Again we have $S_x'\in f^{a_{i}}(U_{\delta',k}^{a_{i}})$
and $\sigma_x'\in U_{\delta',k}^{a_{i}}$. As $a_i\geq i$, by
(\ref{ZikZk2}) we have
$Z^{a_{i}}(j)=Z(j)$ for all $j \leq \nu_i$.
The relation $P(B_\delta^j)\cap P(B_\delta^k)\neq \emptyset$
implies $k\in Z(j)$. Hence
$U_{\delta',j}^{a_i}\cap U_{\delta',k}^{a_i}\neq \emptyset$ and by
Lemma \ref{inclusionnoncompact} b) $U_{\delta',k}^{a_{i}}\subset U_{r,j}^{a_{i}}$.
As $h(x)$ intersects $f^{a_i}(U_{r,j}^{a_i})$ in exactly one point, we
conclude $S_x'=S_x$. Hence $\varphi^{a_i}$ is well-defined. \\ \\
In that way we define for any $i\in \N$ a map $\varphi^{a_{i}}$. Moreover
we may choose $(a_i)_{i\in\N}$ to be
strictly increasing. Passing to the subsequence $f^{a_{i}}$ and simply writing
$f^i$ for it, we obtain maps
$\varphi^i:\bigcup_{j=1}^{\nu_{i}}P(B_\delta^j)\rightarrow
\bigcup_{j=1}^{\nu_{i}} U_{r,j}^i$. \\ \\ \\
\textbf{Step 2:\; \boldmath $\varphi^i(P(B_{\frac{\hat{\delta}}{6}}^j))\subset
U_{\frac{\hat{\delta}}{5},j}^i\subset
\varphi^i(P(B_{\frac{\delta}{3}}^j))\subset
\varphi^i(P(\overline{B_{\frac{\delta}{3}}^j})) \subset
U_{\frac{\delta}{2},j}^i$\unboldmath} \\ \\
Fix $i\in \N$.
As in Lemmas \ref{surj} and \ref{help},
using convergence of graph systems, we may
choose $a_i\in \N$ sufficiently large such that for all $j \leq \nu_i$
\begin{itemize}
\item $\varphi^{a_{i}}(P(B_{\frac{\hat{\delta}}{6}}^j))\subset
U_{\frac{\hat{\delta}}{5},j}^{a_{i}}$,
\item $\varphi^{a_{i}}(P(\overline{B_{\frac{\delta}{3}}^j}))\subset
U_{\frac{\delta}{2},j}^{a_{i}}$,
\item $U_{\frac{\hat{\delta}}{5},j}^{a_{i}}\subset
\varphi^{a_{i}}(P(B_{\frac{\delta}{3}}^j))$. \end{itemize} Doing
so for all $i\in \N$, we may choose $(a_i)_{i\in\N}$ to be strictly
increasing. Denoting the subsequences $f^{a_{i}}$ and $\varphi^{a_{i}}$ simply by $f^i$
and $\varphi^i$, we obtain a sequence with
\beqnn \varphi^i(P(B_{\frac{\hat{\delta}}{6}}^j))\subset
U_{\frac{\hat{\delta}}{5},j}^i\subset
\varphi^i(P(B_{\frac{\delta}{3}}^j))\subset
\varphi^i(P(\overline{B_{\frac{\delta}{3}}^j})) \subset
U_{\frac{\delta}{2},j}^i  \hspace{6mm} \text{ for all } j\leq \nu_i. \eeqnn
\textbf{Step 3:\; Construction of diffeomorphisms\! \boldmath$\phi^i$\unboldmath} \\ \\
Consider the maps
$\varphi^{i}:\bigcup_{j=1}^{\nu_{i}}P(B_\delta^j)\rightarrow
\bigcup_{j=1}^{\nu_{i}} U_{r,j}^i$. \\ \\
By Step 2 we have $U_{\hat{\delta}/5,j}^i\subset\varphi^i(P(B_{\delta/3}^j))$ for any $j\leq\nu_i$.
As $Q^i$ is also a $\frac{\hat{\delta}}{5}$-net for $f^i$, we conclude \beqnn
(f^i)^{-1}(\hat{B}_i)\subset \bigcup_{j=1}^{\nu_i}U_{\hat{\delta}/5,j}^i\,,\eeqnn
so in particular \beqn
(f^i)^{-1}(\hat{B}_i)\subset \varphi^i(\bigcup_{j=1}^{\nu_{i}}P(B_{\delta/3}^j)). \label{deltadiv3ineq} \eeqn
Define
$U^i:=(\varphi^i)^{-1}((f^i)^{-1}(\hat{B}_i))$.
As $f^i$ and $\varphi^i$ are continuous, we conclude that $U^i$ is open.
Restricting $\varphi^i$ to
$U^i$ yields maps \beqnn \phi^i:U^i\rightarrow (f^i)^{-1}(\hat{B}_i).
\eeqnn
By definition $\phi^i$ is surjective. By Step 2 we have
$\varphi^i(P(B_{\delta/3}^j))\subset U_{\delta/2,j}^i$ for $j\leq \nu_i$.
This is all we need to follow
the arguments of Lemma \ref{locinj} and to conclude that
$\varphi^i$ is injective on $P(B_\delta^j)$ for every $j\leq\nu_i$. Moreover, we know by Step 2 that
$\varphi^i(P(B_{\hat{\delta}/6}^j))\subset U_{\hat{\delta}/5,j}^i$ for every $j\leq \nu_i$,
which enables us to follow the arguments of Lemma \ref{inj} in order to conclude that
$\varphi^i$ is injective on all of $\bigcup_{j=1}^{\nu_{i}}P(B_\delta^j)$. Hence also $\phi^i$
is injective on $U^i$.
\\ \\
With the arguments of Lemma \ref{compimmersion} we may show that
$f^i\circ \phi^i$ are immersions and finally that $\phi^i:U^i\rightarrow (f^i)^{-1}(\hat{B}_i)$
are diffeomorphisms. \\ \\ \\
\textbf{Step 4:\; \boldmath$U^i \subset\!\subset U^{i+1}$ and $\bigcup_{i=1}^\infty
U^i=M$\unboldmath} \\ \\
First observe that by Step 2 for $1\leq j \leq \nu_i$, as $\nu_i<\nu_{i+1}$ we have \beqnn
\varphi^{i+1}(x)\in U_{\delta/2,j}^{i+1} \hspace{3mm}\text{ for } x \in
P(\overline{B_{\delta/3}^j}). \eeqnn Hence by
Lemma \ref{distancedeltanet} a) \beqnn
f^{i+1}\circ \varphi^{i+1}(x) \in \hat{B}_{i+1} \hspace{3mm}\text{ for }
x \in \bigcup_{j=1}^{\nu_{i}} P(\overline{B_{\delta/3}^j}). \eeqnn By construction of the
sets $U^i$ this means \beqnn \bigcup_{j=1}^{\nu_{i}}
P(\overline{B_{\delta/3}^j}) \subset U^{i+1}. \eeqnn
Moreover, by (\ref{deltadiv3ineq}) and as $\varphi^i$ is injective on $\bigcup_{j=1}^{\nu_{i}}P(B_\delta^j)$,
we conclude \beqnn
U^i\subset\bigcup_{j=1}^{\nu_i}P(B_{\delta/3}^j), \eeqnn
hence \beqnn
\overline{U^i}\subset
\bigcup_{j=1}^{\nu_{i}}P(\overline{B_{\delta/3}^j}) \subset U^{i+1}.
\eeqnn Moreover we observe \beqnn \bigcup_{i=1}^\infty U^i =
\bigcup_{i=0}^\infty U^{i+1} \supset \bigcup_{i=0}^\infty
\bigcup_{j=1}^{\nu_{i}}P(\overline{B_{\delta/3}^j})\supset \bigcup_{i=1}^\infty
P(B_{\delta/3}^j)=M, \eeqnn hence $\bigcup_{i=1}^\infty U^i=M$. \\ \\
A short \emph{technical remark:} If we pass another time to subsequences $f^{a_i}$,
$\phi^{a_{i}}$ of $f^i$, $\phi^i$ (as will be done in Step 5 and 6), we shall restrict
$\phi^{a_{i}}$ to $V^{a_{i}}:=(\varphi^{a_{i}})^{-1}((f^{a_{i}})^{-1}(\hat{B}_i))\subset U^{a_{i}}$.
Simply writing $f^i,\phi^i,V^i$ instead of $f^{a_{i}},\phi^{a_{i}},V^{a_{i}}$ and after that $U^i$ instead of
$V^i$, we again obtain a sequence of sets $U^i$ with $U^i\subset\!\subset U^{i+1}$ and
$\bigcup_{i=1}^\infty U^i=M$. Conversely we could also first choose all subsequences and define
the sets $U^i$ afterwards.
\\ \\ \\
\noindent \textbf{Step 5:\; A
subsequence with \boldmath
$f^i\circ\phi^i \rightarrow f$ in $C^0$ \unboldmath} \\ \\
For $i\geq k$ we set $\Theta_k^i:=\|f^i\circ
\varphi^i-f\|_{C^{0}(\bigcup_{j=1}^{\nu_{k}}P(B_\delta^j))}$. Now
we fix $k \in \N$. By the convergence argument of the compact case
there is a subsequence $b_i$ with \beqnn \Theta_k^{b_{i}}\rightarrow
0 \text{ as } i\rightarrow\infty. \eeqnn
Hence we may choose a strictly increasing sequence $a:\N\rightarrow \N$
(in particular $a_i\geq i$ for all $i\in\N$) with \beqnn
\Theta_i^{a_{i}}<\frac{1}{i}. \eeqnn
Passing to the subsequence $f^{a_{i}}\circ \varphi^{a_{i}}$ and denoting this sequence simply
by $f^i\circ\phi^i$, we obtain a sequence with
\beqnn
\Theta_i^i=\|f^i\circ
\varphi^i-f\|_{C^{0}(\bigcup_{j=1}^{\nu_{i}}P(B_\delta^j))}
< \frac{1}{i} \rightarrow 0 \text{ as } i\rightarrow\infty. \eeqnn
Restricting $\varphi^i$ to $U^i$ and using the definition of $\phi^i$, we finally obtain
\beqnn \|f^i\circ \phi^i-f\|_{C^{0}(U^i)}\rightarrow 0 \text{ as } i\rightarrow \infty. \eeqnn
\\
\noindent \textbf{Step 6:\; Higher order convergence} \\ \\
We like to find another subsequence such that $f^i\circ \phi^i\rightarrow f$
locally smoothly. This means convergence with respect to
the weak topology $C_{_W}^{\infty}(M,\R^n)$ as defined in \cite{hirsch}, p.\ 34--36,
which in our case is the same as convergence of
$f^i\circ \phi^i\circ \varphi^{-1}$ to $f\circ \varphi^{-1}$ in $C^\K(\varphi(U),\R^n)$
for any chart $(\varphi,U)$ of the atlas $\mathfrak{A}$ constructed above and for any
$\K\in \N_0$. \\ \\
Let $\bar{\rho}$ be a sequence with $\bar{\rho}_j=r_l$ for all $j,l\in \N$ with
$\nu_{l-1}<j\leq \nu_l$. Then for $j\leq \nu_l$, using
$\|\nabla^\K A^i\|_{L^\infty(\hat{B}_{l+1})}\leq C_\K(l+1)$, Lemma \ref{higherderivatives}
implies \beqn \label{bounds1higherder}
\|u_j^i\|_{C^{\K+2}(B_{\bar{\rho}_j})}\,\leq\: C_l^\K, \eeqn
where $C_l^\K$ is a constant depending on $r_1,\, \alpha,\, C_0(l+1),\ldots,\,C_\K(l+1)$\;
(here we use $\bar{\rho}_j\leq r_1$ as $r$ is assumed to be decreasing). \\ \\
Again let $\rho$ be a sequence with $\rho_j=\delta_l$ for all $j,l\in \N$ with
$\nu_{l-1}<j\leq \nu_l$.
Using (\ref{bounds1higherder}) and Theorem \ref{konvergenz}, we may choose
successively subsequences for any $l$ and pass to the diagonal sequence in order to obtain \beqnn
 f^i\circ \phi^i\circ \varphi_j^{-1}\rightarrow f\circ \varphi_j^{-1}
 \hspace{6mm} \text{ in } C^1(B_{\rho_j},\R^n)
\eeqnn
for all $j\in \N$ and charts $\varphi_j:=\varphi_{B_{\rho_{j}}}^j\!\!:\,P(B_\delta^j)\rightarrow B_{\rho_{j}}$. \\ \\
Starting from this sequence, we may choose again successively subsequences for any $\K$ and pass to the diagonal sequence
in order to obtain \beqnn
f^i\circ \phi^i\circ \varphi_j^{-1}\rightarrow f\circ \varphi_j^{-1}
 \hspace{6mm} \text{ in } C^\K(B_{\rho_j},\R^n) \eeqnn
for all $j\in \N$ and all $\K\in \N_0$. \pagebreak \\
\noindent \textbf{Step 7:\; Bounds for the limit} \\ \\
We like to show that \beqnn
\mu(\hat{B}_R) &\leq& C(R) \;\;\text{ for any } R>0, \\
\|\nabla^\K A\|_{L^{\infty}(\hat{B}_{R})} &\leq& C_\K(R) \;\text{ for any } R>0 \text{ and }\,\K \in
\N_0, \eeqnn
where $\mu=f(\mu_g)$ and $A$ is the second fundamental form of $f$. \\ \\
For the first inequality let $R>0$. Let $\varepsilon>0$.
Choose $\tilde{R}$ with $0<\tilde{R}<R$ and \beqn \label{stepbounds1}
\mu(\hat{B}_R)\leq \mu(\hat{B}_{\tilde{R}})+\varepsilon \eeqn
(which is always possible).
By Theorem \ref{compactness1} we have \beqnn
\|f^i\circ \phi^i-f\|_{C^0(U^i)}\rightarrow 0 \hspace{5mm} \text{ as } i\rightarrow\infty, \eeqnn
where $U^i=(f^i\circ \phi^i)^{-1}(\hat{B}_i)$. This yields \beqn \label{stepbounds2}
f^{-1}(\hat{B}_{\tilde{R}})\,\subset\, (f^i\circ \phi^i)^{-1}(\hat{B}_R) \hspace{10mm} \text{ for i sufficiently large.}
\eeqn
Denoting by $\tilde{g}^i$ the metric induced by $f^i\circ \phi^i$, we have \beqn \label{stepbounds3}
\mu_{\tilde{g}^i}((f^i\circ \phi^i)^{-1}(\hat{B}_R))\,=\,\mu_{g^i}((f^i)^{-1}(\hat{B}_R))
\,=\,\mu^i(\hat{B}_R)\,\leq\, C(R). \eeqn
Moreover \beqn \label{stepbounds4}
f^i\circ \phi^i\rightarrow f \hspace{5mm} \text{ on } f^{-1}(\hat{B}_{\tilde{R}})\subset M
\:\text{ in } C^1. \eeqn
Using (\ref{stepbounds4}) in the first line, (\ref{stepbounds2}) in the second,
and (\ref{stepbounds3}) in the third, we obtain \beqnn
\mu(\hat{B}_{\tilde{R}})=\mu_g(f^{-1}(\hat{B}_{\tilde{R}}))
&=& \lim_{i\rightarrow\infty} \mu_{\tilde{g}^i}(f^{-1}(\hat{B}_{\tilde{R}})) \\
&\leq& \limsup_{i\rightarrow\infty} \mu_{\tilde{g}^i}((f^i\circ \phi^i)^{-1}(\hat{B}_R)) \\
&\leq& C(R). \eeqnn
With (\ref{stepbounds1}) this implies $\mu(\hat{B}_R)\leq C(R)+\varepsilon$. As this is true for
any $\varepsilon>0$, we finally conclude $\mu(\hat{B}_R)\leq C(R)$. \\ \\
The bound $\|\nabla^\K A\|_{L^{\infty}(\hat{B}_{R})} \leq C_\K(R)$ is shown in the same way,
using the locally smooth convergence. Note that the first bound would also follow from
Corollary \ref{measureconvergence}, which is shown below.
This completes the proof of Theorem \ref{compactness1}. \hfill $\square$ \\ \\[1mm]
\begin{rem} \label{cooprem}
In \cite{cooper}, the projection for the construction of the diffeomorphisms $\phi^i$ is not carried out.
In the cited paper, one considers sets
$W_{l,i}=\bigcup_{j=1}^{K^l} U^i_{3\delta/4,j}\subset M^i$, where $K^l$ is a constant. Passing to
a subsequence, the corresponding Euclidean isometries $A_j^i$ converge for $1\leq j\leq K^l$.
Choosing $i,i'$ large enough,
the corresponding graph systems of $W_{l,i}$ and $W_{l,i'}$ are close to each other.
It is concluded, that $f^i(W_{l,i})$ is a graph over $f^{i'}(W_{l,i'})$
and that, therefore, $W_{l,i}$ and $W_{l,i'}$ are diffeomorphic.
However, this conclusion is false. In fact, one can easily construct two sets $W_{l,i}$ and $W_{l,i'}$ that are arbitrarily
close in the sense of graph systems, but not diffeomorphic --- for example $W_{l,i}=\mathbb{S}^1$, and $W_{l,i'}$
a spiral that is close to $\mathbb{S}^1$ but diffeomorphic to an open interval. Similarly, there are
counterexamples where both $W_{l,i}$ and $W_{l,i'}$ are noncompact; also assuming a property of the intersections
of the graphs as in Lemma \ref{existencedeltanet} b)
does not suffice for the conclusion. Instead, one has to construct projections between the immersions
in order to obtain diffeomorphisms between \emph{appropriate subsets} of $W_{l,i}$ and $W_{l,i'}$. The same
is needed in order to obtain $C^{\K}$-convergence.
This was done by Langer in
\cite{langer}, and by the author in the present paper.
\end{rem}
\vspace{5mm}
\textbf{Proof of Corollary \ref{measureconvergence}:} \\
The measures $\mu^i$ converge to $\mu$ in $C_c^0(\R^n)'$
if and only if
\beqnn
\lim_{i\rightarrow\infty}\int_{\R^n}f\,d\mu^i=\int_{\R^n}f\,d\mu \hspace{4mm}
\text{ for all } f\in C_c^0(\R^n). \eeqnn
By \cite{evans}, p.\ 54, Theorem 1, this is equivalent to the inequalities \beqnn
\limsup_{i\rightarrow\infty}\mu^i(K)&\leq& \mu(K) \hspace{13mm}\text{ for each compact set }
K\subset\R^n \:\text{ and } \\
\mu(U)&\leq& \liminf_{i\rightarrow \infty}\mu^i(U) \hspace{2.03mm} \text{ for each open set }
U\subset \R^n. \eeqnn
We now will show these two inequalities. \\ \\
So let $K\subset \R^n$ be compact. Let $V\subset \R^n$ be open with $K\subset V$.
By Theorem \ref{compactness1} we have \beqnn
\|f^i\circ \phi^i-f\|_{C^0(U^i)}\rightarrow 0 \hspace{5mm} \text{ as } i\rightarrow\infty, \eeqnn
where $U^i=(f^i\circ \phi^i)^{-1}(\hat{B}_i)$. Hence \beqnn
(f^i\circ \phi^i)^{-1}(K)\subset f^{-1}(V) \hspace{5mm} \text{ for } i \text{ sufficiently large.} \eeqnn
Thus we get, denoting by $\tilde{g}^i$ the metric induced by $f^i\circ \phi^i$, \beqnn
\mu^i(K)\,=\,\mu_{g^i}((f^i)^{-1}(K))\,=\,\mu_{\tilde{g}^i}((f^i\circ \phi^i)^{-1}(K))
\,\leq\, \mu_{\tilde{g}^i}(f^{-1}(V)). \eeqnn
Letting $i\rightarrow \infty$ yields \beqnn
\limsup_{i\rightarrow \infty} \mu^i(K)\leq \mu_g(f^{-1}(V))=\mu(V). \eeqnn
As $\mu(K)=\inf \{\mu(W):\, W\text{ open},\, K\subset W\}$ by Theorem 1.3 in \cite{simon},
we finally obtain
$\limsup_{i\rightarrow \infty}\mu^i(K)\leq \mu(K)$. \\ \\
Next let $U\subset \R^n$ be open. Let $C\subset \R^n$ be compact with $C\subset U$. Then
\beqnn
f^{-1}(C)\subset (f^i\circ \phi^i)^{-1}(U) \hspace{5mm} \text{ for } i \text{ sufficiently large.} \eeqnn
This implies \beqnn
\mu_{\tilde{g}^i}(f^{-1}(C))\leq \mu_{\tilde{g}^i}((f^i\circ \phi^i)^{-1}(U))
=\mu_{g^i}((f^i)^{-1}(U))=\mu^i(U). \eeqnn
Again letting $i\rightarrow \infty$ yields \beqnn
\mu(C)=\mu_g(f^{-1}(C))\leq \liminf_{i\rightarrow \infty} \mu^i(U). \eeqnn
As $\mu(U)=\sup\{\mu(E):\, E\text{ compact},\, E\subset U\}$ by Remark 1.4 in \cite{simon}, we obtain
$\mu(U)\leq \liminf_{i\rightarrow \infty}\mu^i(U)$, which proves Corollary \ref{measureconvergence}.
\hfill $\square$
\vspace{1cm} \\
Finally we would like to give some generalizations of Theorem \ref{compactness1}. First we remark that the assumptions
(\ref{comp1eq1}) and (\ref{comp1eq2}) can be weakened as follows: Let $f^i:M^i\rightarrow \R^n$ be as in Theorem \ref{compactness1}
with $f^i(M^i)\cap K\neq \emptyset$ for a compact set $K\subset \R^n$.
Let $(R_i)_{i\in \N}$ be a sequence in $\R_{>0}$ with $R_i\rightarrow \infty$ as $i\rightarrow \infty$,
and assume \begin{align}
\hspace{4.8cm} \mu^i(B_R) \:&\leq\: C(R) \;\;\hspace{0.4mm} \text{ for any } R<R_i, \hspace{4cm} &(i) \nonumber\\
\|\nabla^\K A^i\|_{L^{\infty}(B_{R})} \:&\leq \:C_\K(R) \;\text{ for any } R<R_i \text{ and }\,\K \in
\N_0. &(ii) \nonumber \end{align}
Then the same statement as is Theorem \ref{compactness1} holds. The bounds (\emph{i}) and (\emph{ii}) for the
limit, that is $\mu(B_R) \:\leq\: C(R)$ and $\|\nabla^\K A\|_{L^{\infty}(B_{R})} \:\leq \:C_\K(R)$,
hold for any $R>0$. This statement is needed in \cite{link}. The proof is essentially the same as for Theorem \ref{compactness1}. \\ \\
A further reaching generalization is to consider proper immersions into open subsets $\Omega\subset \R^n$:
\begin{cor} \label{coropen} Let $f^i:M^i\rightarrow \Omega$ be a sequence of proper immersions, where $M^i$
is an $m$-manifold without boundary, $\Omega\subset \R^n$ open, and $f^i(M^i)\cap \mathcal{C}\neq \emptyset$ for a compact set
$\mathcal{C}\subset \Omega$. Assume \begin{align}
\hspace{4.8cm} \mu^i(K) \:&\leq\: C(K) \;\;\hspace{0.4mm} \text{ for any } K\subset\Omega \text{ compact}, \hspace{2.65cm} &(i\hspace{0.5mm}') \nonumber\\
\|\nabla^\K A^i\|_{L^{\infty}(K)} \:&\leq \:C_\K(K) \;\text{ for any } K\subset\Omega \text{ compact and }\,\K \in
\N_0. &(ii\hspace{0.5mm}') \nonumber \end{align}
Then there exists a proper immersion $f:M\rightarrow \Omega$, where $M$ is again an $m$-manifold without boundary, such
that after passing to a subsequence there are
diffeomorphisms \beqnn \phi^i:U^i\rightarrow (f^i)^{-1}(\Omega^i)\subset M^i, \eeqnn
where $\Omega^i\subset\Omega$, $U^i\subset M$ are open sets with $\Omega^i\subset\!\subset \Omega^{i+1}$,
$U^i\subset\!\subset U^{i+1}$
and $\Omega=\bigcup_{i=1}^{\infty}\Omega^i$, $M=\bigcup_{i=1}^{\infty}U^i$, such that
$\|f^i\circ \phi^i-f\|_{C^{0}(U^i)}\rightarrow 0$,
and moreover $f^i\circ \phi^i\rightarrow f$ locally smoothly on $M$. \\ \\
Moreover, the immersion $f$ also satisfies (i\hspace{0.5mm}$'$\!) and (ii\hspace{0.5mm}$'$\!), that is
$\mu(K) \leq C(K)$ and $\|\nabla^\K A\|_{L^{\infty}(K)}\linebreak \leq C_\K(K)$.
\end{cor}
\textbf{Proof:} \\
We set $V^i:=B_i(0)\cap \Omega_{1/i}$\,, where $\Omega_\delta:=\{x\in \Omega: \text{dist}(x,\partial \Omega)>\delta\}$.
As $\mathcal{C}$ is compact, there is an $i_0\in \N$ with $\mathcal{C}\subset V^i$ for all $i\geq i_0$. We set $\Omega^i:=V^{i_0+i}$.
As $\overline{\Omega^i}\subset \Omega$ is compact, on $\Omega^i$ we have uniform bounds for the volume and the second fundamental form.
Now we can proceed as in the proof of Theorem \ref{compactness1}: For a given immersion $f:M\rightarrow\Omega$ and $p\in M$ define
$\bar{p}=\bar{p}(f):=\min\{j\in \N:f(p)\in \Omega^j\}$. Now
always replace the balls $\hat{B}_i$ by the sets
$\Omega^i$ (for example in Definition \ref{decomp2}). Following step by step the arguments of
Theorem \ref{compactness1}, the corollary follows. \hfill $\square$
\\ \\
Although in the corollary above the target $\Omega$ is not necessarily complete, the limit does not only lie in
$\overline{\Omega}$, but even in $\Omega$. This is possible as here we only desire local convergence.
Theorem \ref{compactness1} and Corollary \ref{coropen} can be generalized further. First one could formulate a version
of Corollary \ref{coropen} as in the paragraph preceding this statement. Moreover there are versions for proper
$C^\K$-immersions with \beqnn
\|\nabla^l A^i\|_{L^{\infty}(B_{R})} \leq C_l(R) \;\text{ for any } R>0 \text{ and }\, l
\text{ with } 0\leq l \leq \K-2 \eeqnn
(or with uniform bounds on compact sets $K\subset\Omega$) with convergence in $C^{\K-1}$. Finally
also for the noncompact case there are versions for proper immersions into Riemannian manifolds $N$. For that we again use
an isometric embedding $N\hookrightarrow \R^\nu$. \\ \\
We would like to give an example how Corollary \ref{coropen} can be used. Let us consider the graph situation
where $g^i:\R^m\rightarrow \R^n$, $g^i(x)=(x,u^i(x))$ with $u^i:\R^m\rightarrow \R^k$.
Now assume that there is exactly one point where curvature
concentrates, say in $0\in \R^m$. We cut out a ball $\overline{B_\varrho(0)}$, and obtain proper immersions
$f^i:\R^m\backslash \overline{B_\varrho(0)}\rightarrow \Omega$, $f^i:=g^i|_{\R^m\backslash \overline{B_\varrho(0)}}$ with
$\Omega:=(\R^m\backslash \overline{B_\varrho(0)})\times \R^k$. Note that $f^i$ is \emph{not} proper as a mapping into $\R^n$.
Let us assume that $f^i$ admits uniform bounds on $\nabla^{\K} A^i$. Then we are in a situation to apply Corollary
\ref{coropen}, and to conclude convergence of a subsequence.
Similarly, we can also apply the corollary in the case of graphs defined on annuli $B_R\backslash B_r$.
\end{section}

\end{document}